\documentclass[letterpaper, 10 pt, conference]{ieeeconf}  

\IEEEoverridecommandlockouts                              

\usepackage{enumerate}
\usepackage{lingmacros}
\usepackage{tree-dvips}
\usepackage{amsmath}

\usepackage{amsthm}
\usepackage{graphicx} 
\usepackage{fancybox}
\usepackage{array}
\usepackage{mdwmath}
\usepackage{mdwtab}
\usepackage{setspace}
\usepackage{fixltx2e}
\usepackage{stfloats}
\usepackage{ifdraft}
\usepackage{graphicx} 

\usepackage{url}
\usepackage{t1enc}
\usepackage{float}
\usepackage{lingmacros}
\usepackage{tree-dvips}
\usepackage{wrapfig}
\usepackage{amssymb}
\usepackage{bm}
\usepackage{algorithm} 
\usepackage{algorithmic}
   
\usepackage{epstopdf}
\usepackage{alphalph}
\usepackage{bm}
\usepackage{mathtools}
\usepackage{caption}
\usepackage{subcaption}

\newtheorem{theorem}{Theorem}[section]

\newtheorem{proposition}[theorem]{Proposition}

\title{\LARGE \bf
A UKF-PF based Hybrid Estimation Scheme for Space Object Tracking
}

\author{Dilshad Raihan A.V$^{1}$ and Suman Chakravorty$^{2}$
\thanks{}
\thanks{$^{1}$Dilshad Raihan A.V is a graduate student researcher at the Department of Aerospace Engineering, Texas A\&M University, College Station
        }%
\thanks{$^{2}$Suman Chakravorty is Associate Professor with the Department of Aerospace Engineering, Texas A\&M Univeristy,
        College Station
        }%
}

\begin{document}

\maketitle
\thispagestyle{empty}
\pagestyle{empty}

\begin{abstract}

In this paper, we present a UKF-PF based hybrid nonlinear filter for space object tracking.  Estimating the state and its associated uncertainty, also known as filtering is paramount to the tracking process. The periodicity of the Keplerian orbits and the availability of accurate orbital perturbation models present special advantages in filter design. The proposed nonlinear filter employs an unscented Kalman filter (UKF) estimate the state of the system while measurements are available. In the absence of measurements, the state pdf is updated via a sequential Monte Carlo method. It is demonstrated that the hybrid filter offers fast and accurate performance regardless of orbital parameters used and the amount of uncertainty involved. The performance of the filter under is found to depend upon the number of measurements recorded when the object is within the field of view (FOV) of the sensors.

\end{abstract}

\section{INTRODUCTION}

The number of objects that populate the earth's sky has gone through a great surge over the years. Collision with debris and decommissioned satellites pose a real risk hazarding the feasibility of space operations and satellites\cite{Dona}. This has given rise to a great appeal for the development of accurate estimation schemes for space object tracking. The optimal linear estimator known as the Kalman filter proposed by Rudolf Kalman set the frame work for recursive estimation of  uncertain dynamical systems\cite{Kal},\cite{Bucy}. The Kalman filter furnishes the unbiased  minimum variance estimates when the dynamical system is linear and the uncertainties involved are Gaussian. The perturbed orbital dynamics of the space objects, like many other dynamical systems in nature, is essentially nonlinear. The development of the extended Kalman filter (EKF) set forth attempts  to derive the optimal filter for nonlinear dynamical systems\cite{Rist},\cite{Dau},\cite{Smit}. The EKF involved the linearization of the state transition equations and the observation model at the current estimated state. The errors accumulated due to  linearization and restrictive assumptions it enforced on the nature of uncertainties were major shortcomings of EKF. Improved results were produced when a second order EKF that included second order Taylor series terms was employed\cite{Atha},\cite{Geven}. The emergence of sigma point Kalman filters, specifically the unscented Kalman fitler gave rise to an alternative that doesn't rely on linearization of dynamics\cite{Jul},\cite{Julier},\cite{Wan}. The UKF approximated the state pdf with a set of points, also called as the sigma points, carefully chosen so that the first two moments of the original pdf are captured. It was found to perform better than the EKF. However, the UKF, like the EKF, is a finite dimensional filter that estimates only the first two moments of the pdf which could, in general, require an infinite number of parameters for a full description. Hence, the estimation results of both these filters are suboptimal and can diverge. Handling the non-Gaussianity of the state pdf is especially relevant in problems such as space object tracking wherein no measurements may be registered for extended periods of time.\\

The particle filter (PF) or the sequential Monte Carlo estimator is  a nonlinear filter that doesn't
enforce restrictive assumptions on the nature of pdf or dynamics of the system\cite{Aru}\cite{Gor}. The PF employs a suitably large number of particles constituting a representative ensemble of the state pdf for state estimation. It has been proved that the variance of the particle weights is a nondecreasing function of time\cite{Docet}. While the absence of restrictive assumptions make it a more robust estimator suitable for a general nonlinear filtering problem, particle filters suffer from the curse of dimensionality, i.e. as the dimension of the state space increases, it is bound to fail due to weight depletion unless the number of points are increased exponentially\cite{Beng}. Hence, it becomes computationally expensive as the dimensionality of the problem becomes large\cite{Hua}. Estimation of nonlinear dynamical systems by employing a Gaussian mixture model to approximate the state pdf has also been proposed\cite{Alspach},\cite{Sorenson}. Methods to improve the estimation accuracy of Gaussian sum filters with between-measurement weight updates have also gained much attention recently\cite{Terejanu},\cite{Mars}.\\

In this paper, we propose a novel estimation scheme that combines the advantages of the UKF and PF to produce a fast and accurate nonlinear filter that can be employed for space object tracking. A UKF is used to estimate the state of the object when measurements are available, i.e., during its flight within the observer's field of view. As the object moves out of the FOV, an ensemble of particles are sampled and propagated in time based on the available model of the orbital dynamics. The proposed filter makes the best use of the unique features of the space object tracking problem namely periodicity, minimal process noise and relatively small velocity uncertainty. The new filter is demonstrated to be capable of producing fast and accurate estimates  irrespective of the orbital parameters involved and the level of uncertainty. \\

The remainder of this paper is organized as follows. In section II, the dynamics of a space object is briefly reviewed. Section III discusses the details that pertain to the uncertainty propagation in the orbital dynamics which are relevant to the filter design. A detailed account of the filter design process is provided in section IV. Section V then discusses the results obtained when the proposed scheme was employed in two test cases of space objects in low earth orbits subject to atmospheric drag and $J_{2}$ perturbation. The derivation of a particular result concerning the absence of particle depletion in periodic dynamical systems is presented in the appendix.

\section{Dynamics \&  Measurement models}
In this section, the perturbed dynamics of the orbiting objects is briefly discussed, followed by a description of an angles only measurement model employed to aid state estimation.
\subsection{Dynamics of space objects}
The  acceleration experienced by an object in the inverse square  gravitational field of earth is given by
\begin{equation}
a_{g}=-\frac{GMm\vec{r}}{r^{3}}. \label{earthfield}
\end{equation}
Here G is the universal gravitational constant, r the vector joining the center of earth to the CM of the object and M the mass of earth. The gravitational acceleration as given in eqn~\eqref{earthfield} assumes that the central body is  spherically symmetrical.In reality, earth has a non symmetrical mass distribution with more mass distributed along the equator and and is considered something akin to an oblate ellipsoid. To account for the non-sphericity, the gravitational potential is expanded into a series of spherical harmonics.The dominant perturbation term in the resulting expansion is called the $J_{2}$ harmonic. The perturbing acceleration arising from the $J_{2}$ term, $a_{j_{2}}$ is given by
\begin{equation}
a_{J_{2}}=-\frac{3}{2}J_{2}\frac{\mu }{r^{2}}(\frac{r_{eq}}{r})^{2}\begin{pmatrix}(1-5(\frac{z}{r})^{2})\frac{x}{r}\\(1-5(\frac{z}{r})^{2})\frac{y}{r}\\(3-5(\frac{z}{r})^{2})\frac{z}{r}\end{pmatrix}.
\end{equation}
where x , y , z are the Cartesian coordinates of the CM of the object measured from the centre of earth\cite{Junk}.In addition to this, the orbital motion is also affected by the non-conservative atmospheric drag which may be significant in low earth orbits. The acceleration due to drag force is given by
\begin{equation}
a_{D}=-(\frac{A}{m})C_{d}\rho \frac{v^{2}}{2}\vec{i_{v}}.
\end{equation}
Here m is the mass of the object, A its cross sectional area, $C_{d}$ the drag coefficient, $\rho$ the density, and v is the relative velocity between the atmosphere and orbiting object.
 A simple exponential model may be employed to describe the variation of atmospheric density with altitude according to which
\begin{equation}
\rho (r)=\rho _{0}\exp(-(r-r_{0})/H).
\end{equation}
Here $\rho _{0}$ and $r_{0}$ are reference density and radius. The variable H, known as scale height, is the vertical distance over which the density of the atmosphere reduces by a factor of mathematical constant $e$\cite{Meteo}.
\subsection{Measurement Model}
\begin{figure}[h]
\includegraphics[width=0.5\textwidth,height=0.2\textheight]{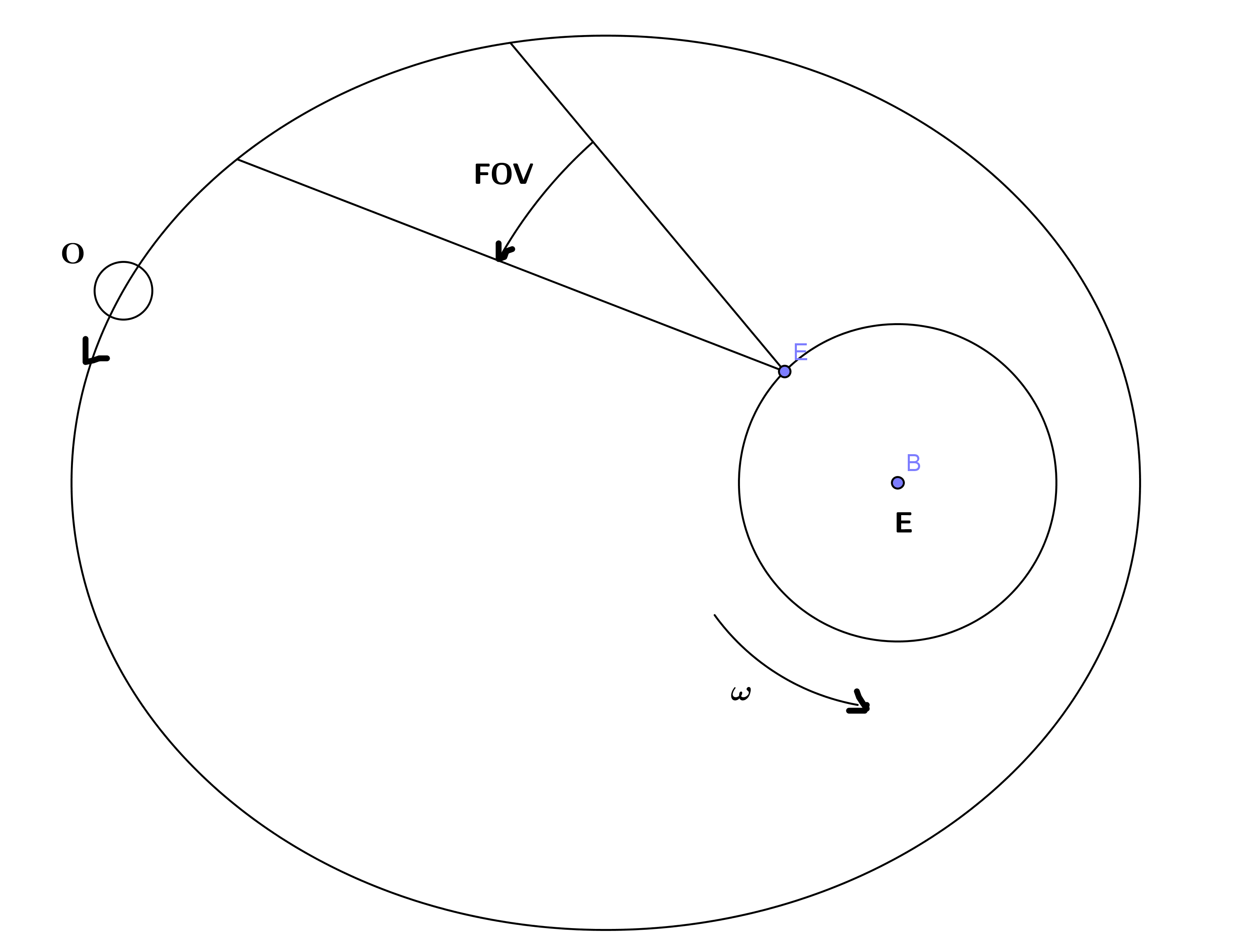}  
\caption{The sensor is fixed on the ground station which defines a non inertial frame that spins with the earth}
\label{fig:spaceobject}
\end{figure}

Let $\vec{r}$ and $\vec{r_{s}}$ be the inertial position vectors of the space object O and the ground station respectively. Then the relative position of the object w.r.t the ground station is given by 
\begin{equation}
\vec{\rho_{i}}= \vec{r}-\vec{r_{s}}.
\end{equation}
The sensor measures the topocentric inclination($\theta$) and right ascension ($\phi$) from a ground station assumed to be located on earth's equator.
 The coordinatization of the relative position vector of the space object with respect to ground station in the station from can be computed by multiplying the inertial vector $\vec{\rho_{i}}$ with the appropriate orthonormal transformation matrix.If the effects due to precession ,nutation etc of the earth are neglected, it would be straightforward to see that ,with respect to the inertial frame,the ground station  is in an elemental rotation about the polar axis.
At t=0, both ground frame and inertial frame are aligned. Assuming a constant spin rate $\omega$ for the earth, the transformation matrix for the ground station at time t can be calculated as
\begin{equation}
C(t)=\begin{bmatrix}
cos\omega t & sin\omega t & 0\\
-sin\omega t & cos \omega t & 0\\
0 & 0& 1
\end{bmatrix}.
\end{equation}If [$\rho_{x}$ $\rho_{y}$ $\rho_{z}$] are the Cartesian coordinates of the object in the ground frame , then the inclination($\theta$) and right ascension($\phi$) are calculated as 
\begin{eqnarray}
\theta=sin^{-1}(\rho_{z} / \rho)\\
\phi=tan^{-1}(\rho_{y} / \rho_{x})
\end{eqnarray}
where $\rho$ equals $\sqrt{\rho_{x}^{2}+\rho_{y}^{2}+\rho_{z}^{2}}$.\\
A zero mean Gaussian measurement noise is assumed with 3.9 arc sec standard deviation in angle measurements. The field of view of the ground station is  limited by 75 degree on the either side in the azimuthal direction and by 90 degree on either side in the polar direction. An illustration of the space object-ground station system is presented in fig~\ref{fig:spaceobject}. Once the space object is inside the FOV, the sensor would attempt to scan the visible part of sky and register a measurement based on a detection probability set at 0.9.

\section{Uncertainty propagation in Orbital dynamics}
Let $X(t,X_{0})$ be the state of a continuous time dynamical system governed by 
\begin{equation}
\dot{X}=f(X)+w(t)
\label{eqn:gover}
\end{equation}
 The stochastic term $w(t)$ represents the effect of modelling uncertainties and external noise. Let $p(0,X_{0})$ be the probability density function describing the initial state. Let $Z_{k}$ be a measurement vector recorded by observing the system. The estimation process may be enhanced by recording a measurement vector, $Z_{k}$, of the system defined by 
 \begin{equation}
Z_{k}=H(X_{k})+\nu _{k}.
\label{eqn:measur}
\end{equation}
 Due to the stochastic terms that appear in equations \ref{eqn:gover} and \ref{eqn:measur}, the knowledge regarding the current configuration of the system can only be presented as a probability density function $p(t,X)$. The primary objective in a filtering problem is to compute the pdf $p(t,X)$ which characterizes the uncertainty involved in state estimates at all times. The time evolution of the pdf of a system subject to random forces is described by a deterministic linear partial differential equation known as The Fokker Planck Kolmogorov (FPK) equation\cite{Risk}. Attractive as it may seem, the FPK equation doesn't permit analytical solutions in most cases. Alternatively, a sampling based approach may be employed wherein the pdf at any time can be approximated by an ensemble of particles sampled from the state space\cite{Douc}.\\
 
\begin{figure}[b]
\begin{flushleft}
\centering
\begin{subfigure}[b]{0.5\textwidth}
\includegraphics[width=\textwidth,height=0.2\textheight]{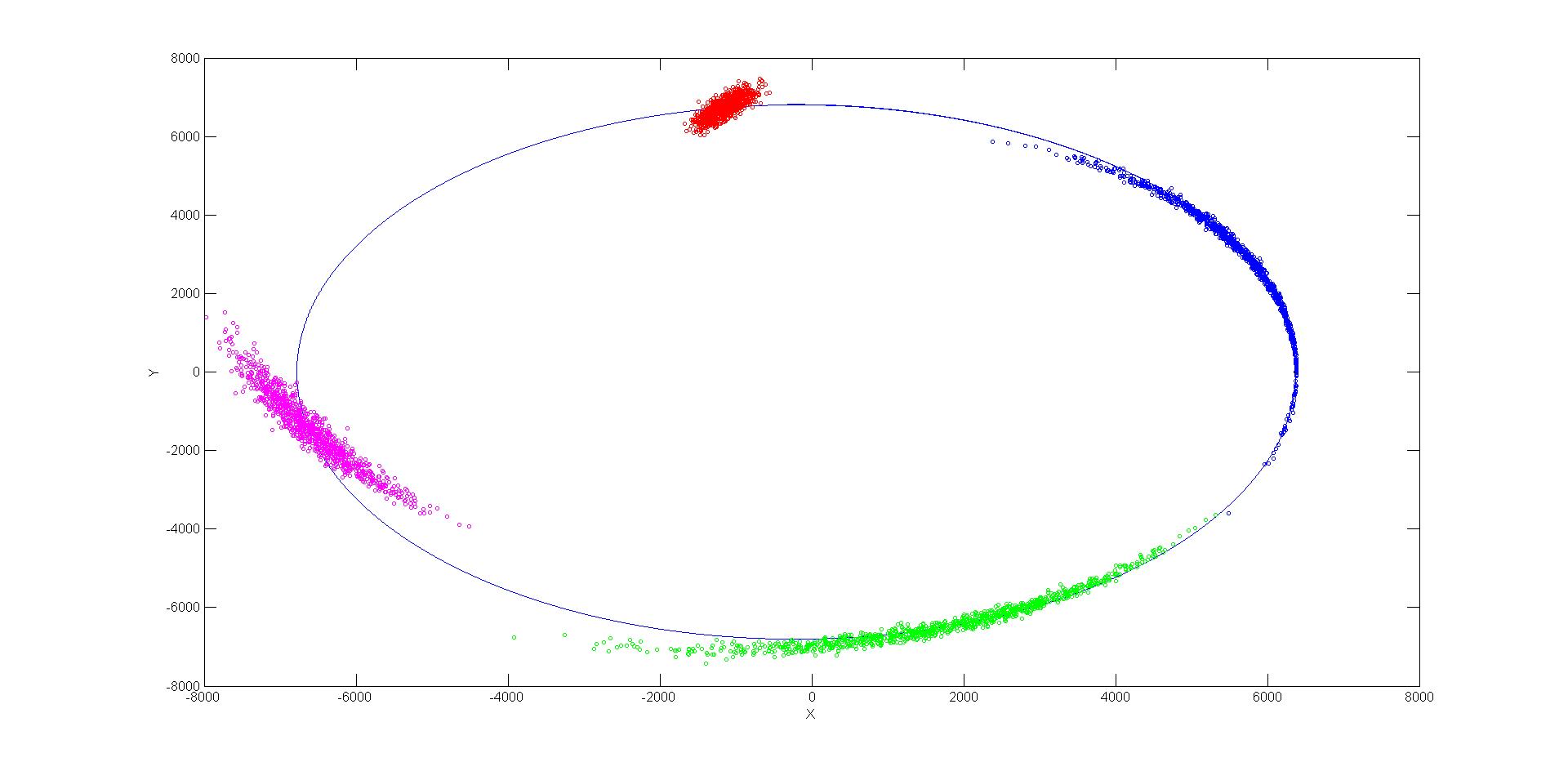}
\caption{} 
\label{fig:ensemble1}
\end{subfigure}
\caption{(a) The ensemble of points represent the marginal distribution of position at times 0, T/4, T/2, 3T/4 and T respectively. (b) The ensemble of particles represent the marginal distribution of velocity at times 0, T/4, T/2, 3T/4, and T respectively}
\label{fig:ensemble}
\end{flushleft}
\end{figure}

To start with, a finite number of particles is selected as a representative ensemble of the initial pdf. This ensemble is then propagated in time based on the equations governing the evolution of the dynamical system. The distribution of states that are occupied by the ensemble at any time is assumed to represent the uncertainty involved in state estimate at that instant. An ensemble of particles, sampled from an initial Gaussian distribution and propagated in time based on the equations of a space object is plotted in fig~\ref{fig:ensemble}. Clearly, regions that are denser with particles represent more probable states. At t=0, the state of the object is assumed to be normally distributed with mean $[6600cos\pi /12 \hspace{6pt} 0 \hspace{6pt} 6600sin\pi/12 \hspace{6pt} 0 \hspace{6pt} 7.8848\hspace{6pt} 0]^{T}$. The  distribution of the particles at  various instants are plotted by propagating the ensemble through the dynamics of the system. The outer boundary of the initial sample resemble an ellipse since the level sets of any Gaussian pdf is elliptical.\\

 As time progresses the initial Gaussian pdf is subjected to a series of nonlinear transformations and the state pdf does not remain Gaussian. The stretched and distorted ensembles validate this conclusion.\\

\indent The sampling based approach doesn't provide a quantitative measure of probability corresponding to any given realization of state vector. That requires a functional form such as a Gaussian mixture model (GMM) of the underlying pdf which may be retrieved from the ensemble of states. A clustering scheme such as Figueiredo-Jain (F-J) algorithm may be employed to arrive at the GMM from the particle sample\cite{Figjain}\cite{Paalanen}. The F-J algorithm uses the minimum message length (MML) criterion to obtain the optimal number of Gaussian components, their weights and corresponding means and covariances, given the ensemble of states. The particle samples plotted in fig~\ref{fig:ensemble1} are clustered with F-J algorithm and the resulting GMM parameters are listed in  Table~\ref{tab:highnoise}. The initial standard deviation in  position and velocity variables are 1 km and 1 km/s respectively in each direction.The number of Gaussian mixture components increases upto 5 with time, evidently due to the growth in uncertainty. The remarkable surge in the trace of the corresponding covariance matrices quantifies the magnitude of this upswing in uncertainty. A graphical representation of the marginal distribution of the state pdf in x and y coordinates is 
\begin{table}[H]
 \centering
 \caption{}
\label{tab:highnoise}
  \begin{tabular}{ |c|c|l| }
  \hline
  \multicolumn{3}{|c|}{High noise} \\
  \hline
  time & modes & \multicolumn{1}{c|}{trace($\times 10^{4}$)} \\
\hline
  0 & 1 & $3.123\times 10^{-4}$  \\
 $1500$ & $4$ & $2.426  $ \space{} $2.210$ \space{} $3.912 $ \space{} $3.155 $\\
   $3000$ & $5$ & $7.807 $ \space{} $4.939$ \space{} $5.195 $ \space{} $28.85$ \space{} $16.23$\\
$4500$ & $4$ & $65.57 $ \space{} $60.30 $ \space{} $57.89 $ \space{} $70.43$\\
$6000$ & $4$ & $9.792 $ \space{} $10.64$ \space{} $14.85$ \space{} $92.70$\\  \hline
\end{tabular} 
\end{table}

presented in figure~\ref{fig:amocodia} \\
 
%

\begin{figure}[H]
       \centering
        \begin{subfigure}[b]{0.2\textwidth}
                \includegraphics[width=\textwidth]{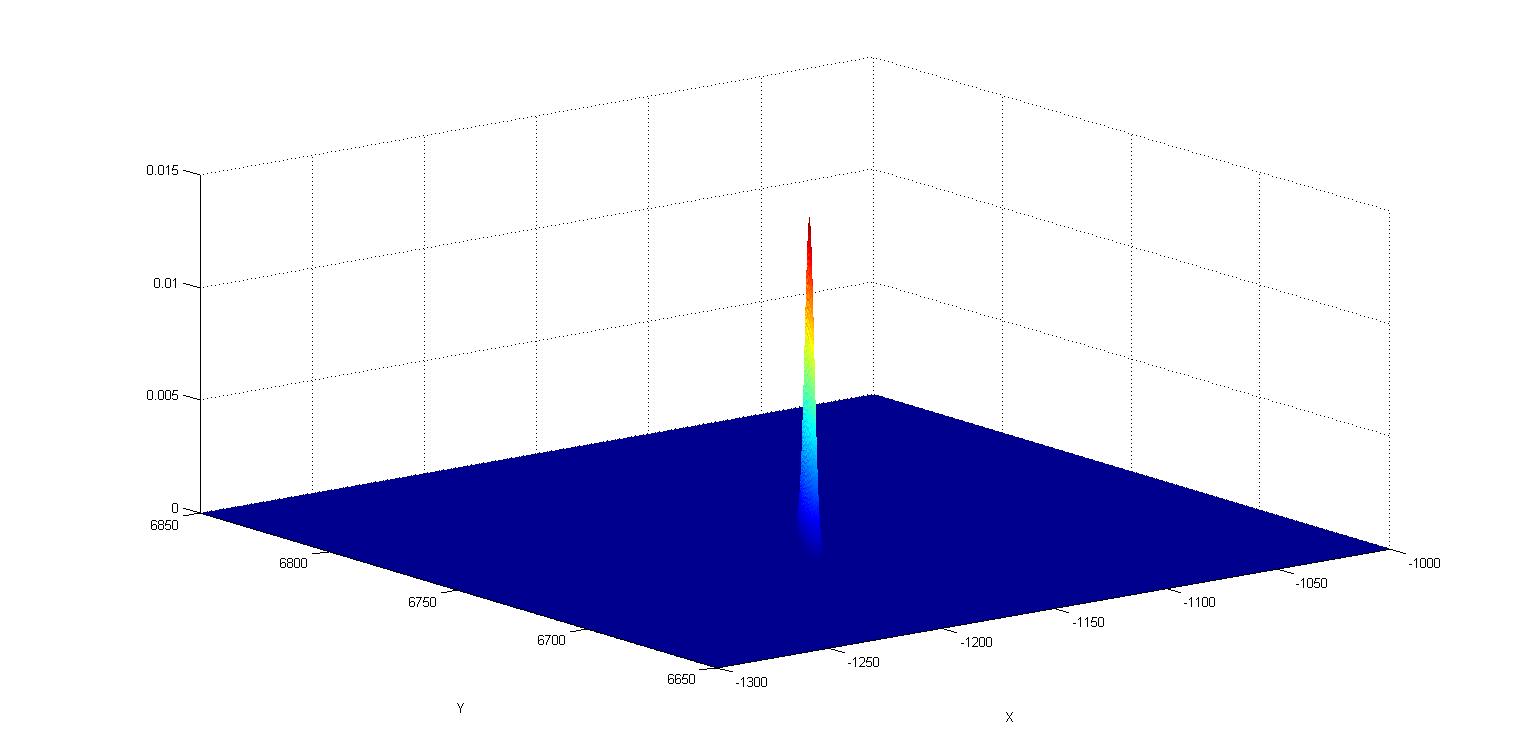}
                \caption{$t=T/4$}
                \label{fig:penddisp}
        \end{subfigure}
        \begin{subfigure}[b]{0.2\textwidth}
                \includegraphics[width=\textwidth]{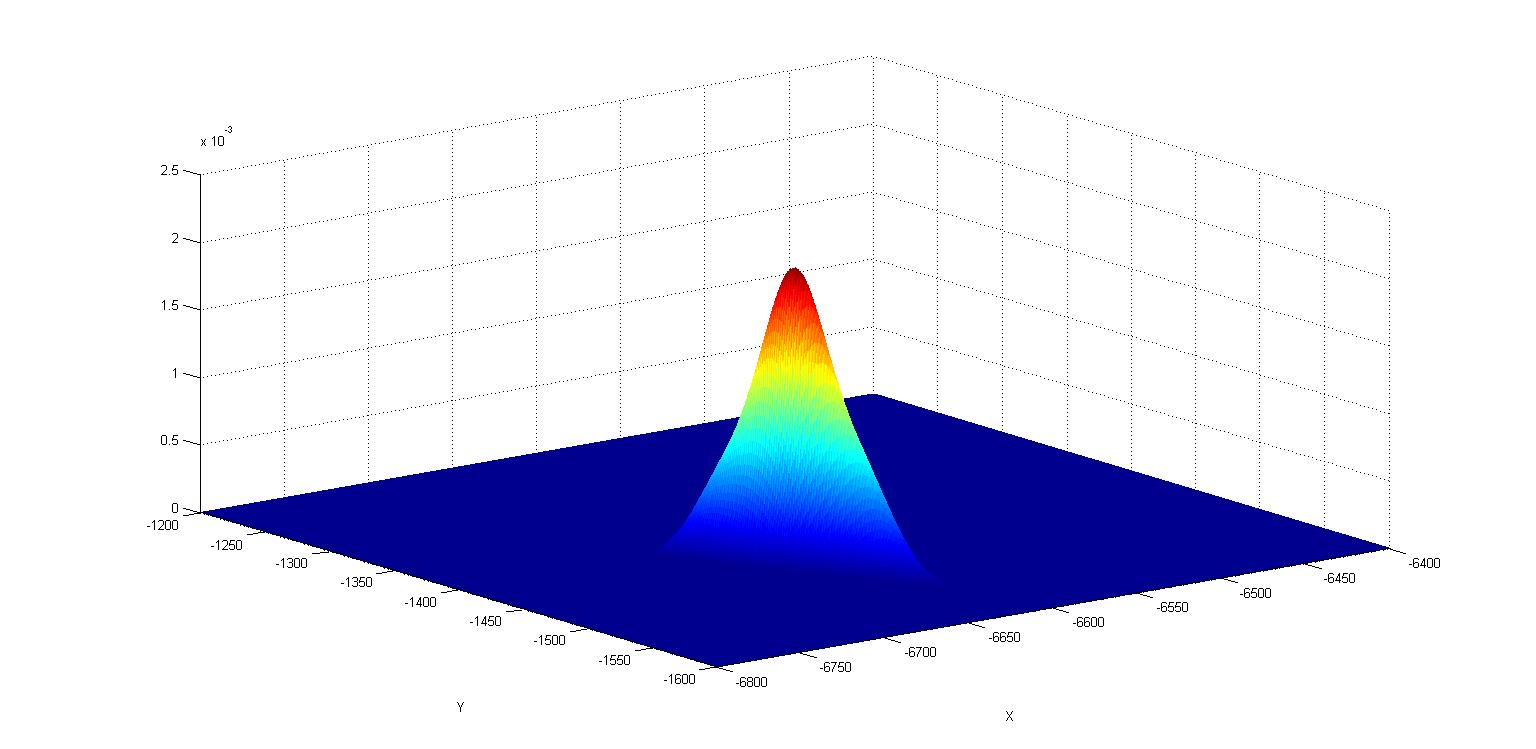}
                                \caption{at  $t=T/2$}
                \label{fig:pendangdisp}
        \end{subfigure}
        
        \begin{subfigure}[b]{0.2\textwidth}
                \includegraphics[width=\textwidth]{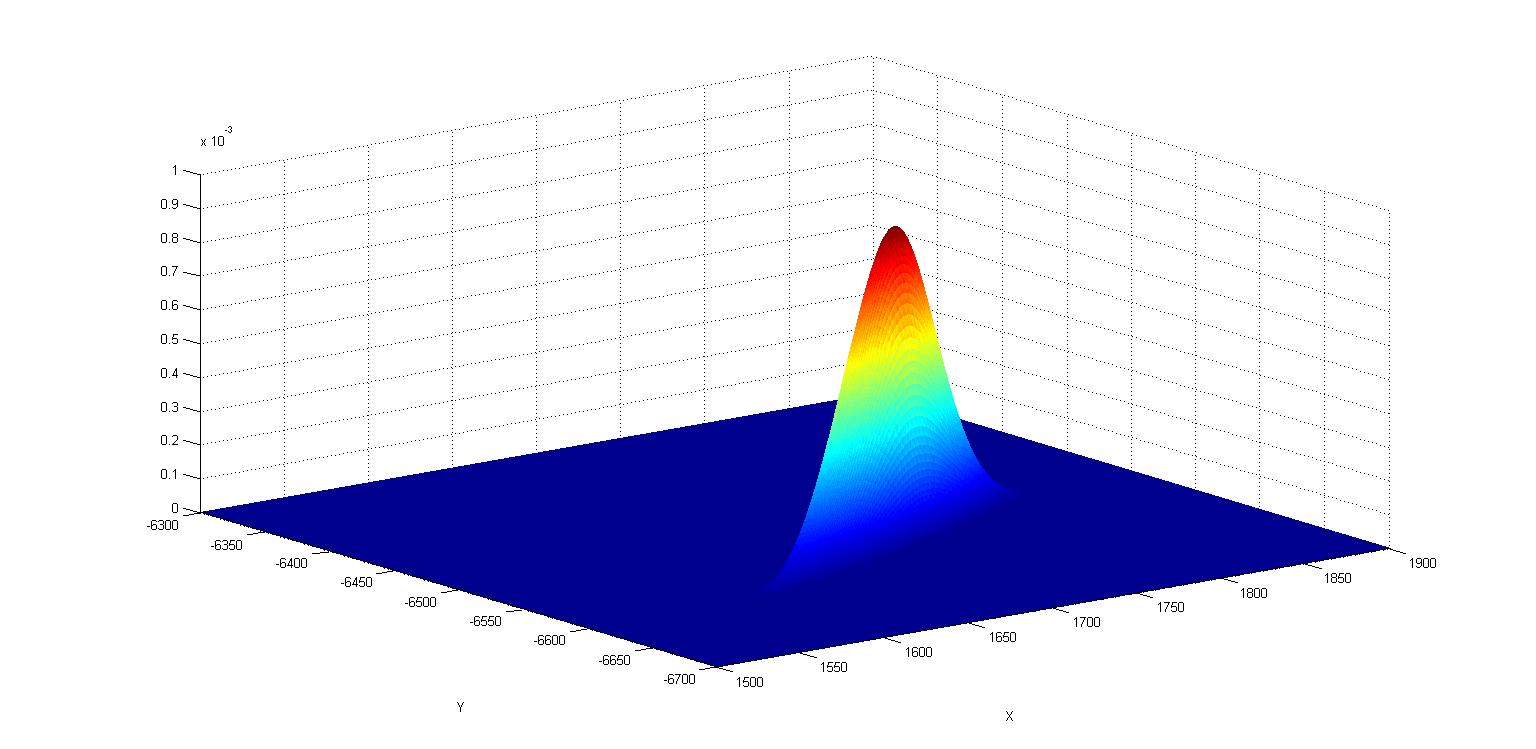}
                \caption{at $t=3T/4$}
                \label{fig:pendvel}
                \end{subfigure}
         \begin{subfigure}[b]{0.2\textwidth}
        \includegraphics[width=\textwidth]{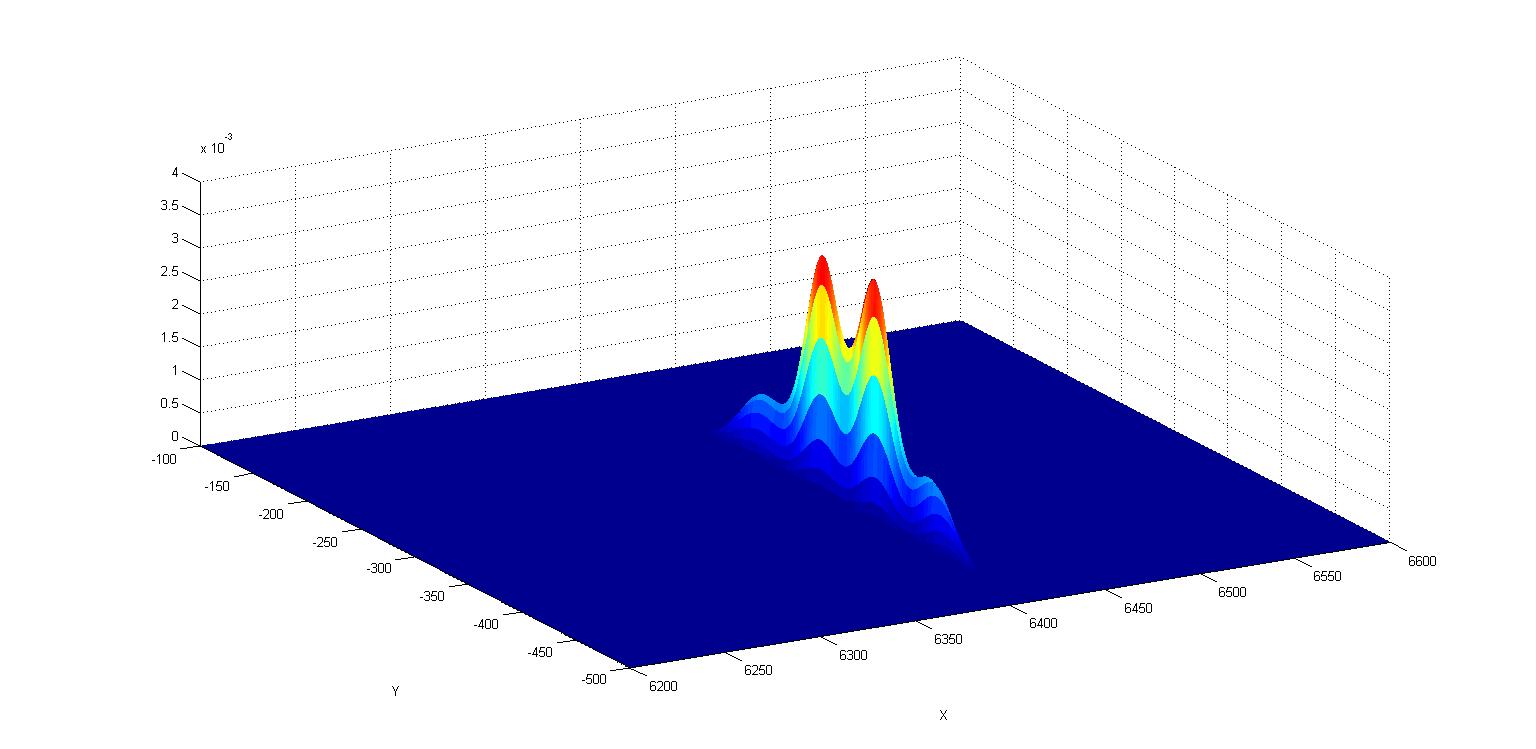}
        \caption{at $t=T$}
        \label{fig:pendangvel}
        \end{subfigure}
              \caption{GMM representation for the ensemble of particles for $t>0$ with lower uncertainty in velocity}
        \label{fig:margdist}
\end{figure}

  Studying the propagation of uncertainty in orbital dynamics reveals three salient features that could be instrumental to the design of an efficient and robust filter. They are:\\
\begin{enumerate}
\item Sensitivity with respect to uncertainty in velocity:  The growth in uncertainty with time is remarkably sensitive towards error in velocity. An ensemble of particles were sampled from a Gaussian distribution keeping
\begin{figure}[H]
       \centering
        \begin{subfigure}[b]{0.2\textwidth}
                \includegraphics[width=\textwidth]{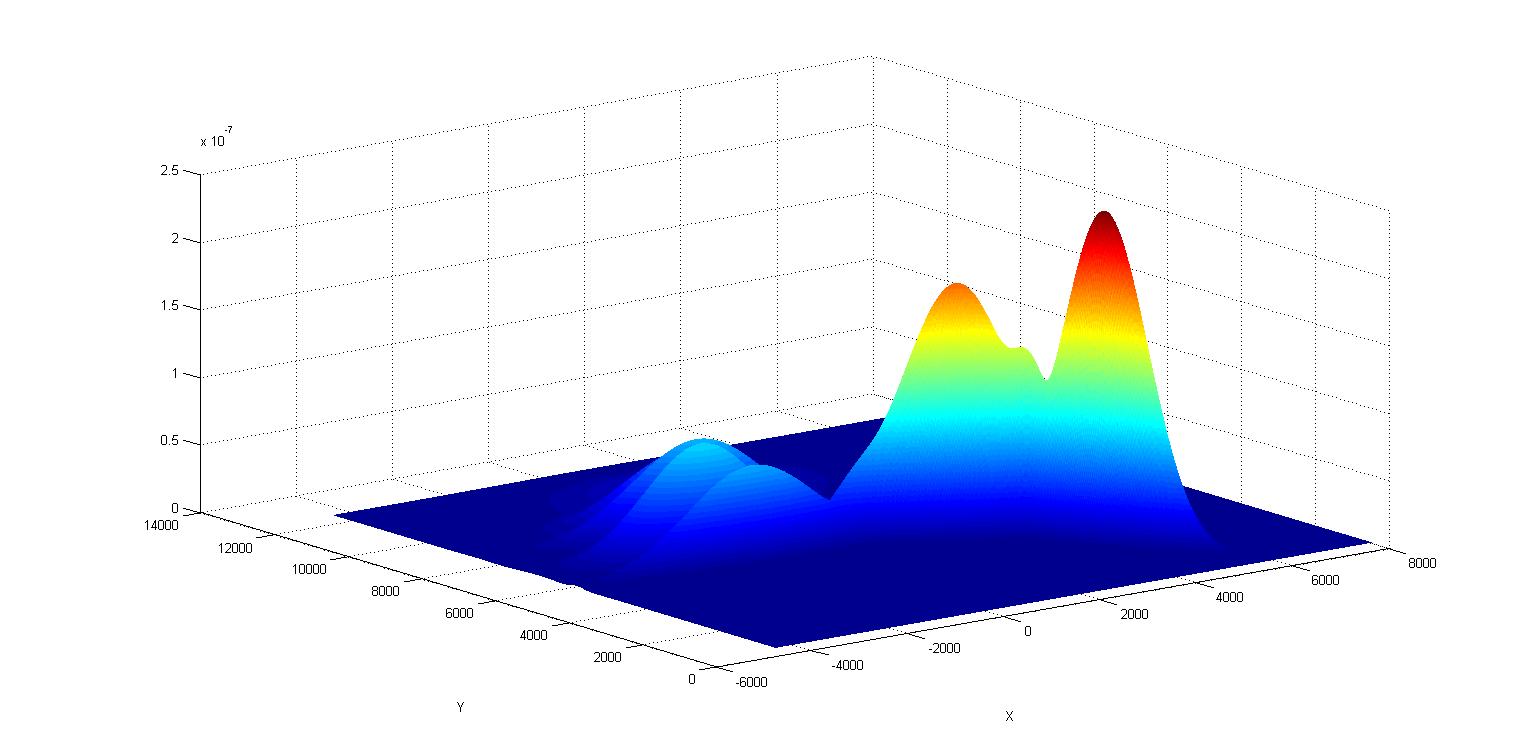}
                \caption{$t=T/4$}
                \label{fig:penddisp}
        \end{subfigure}
        \begin{subfigure}[b]{0.2\textwidth}
                \includegraphics[width=\textwidth]{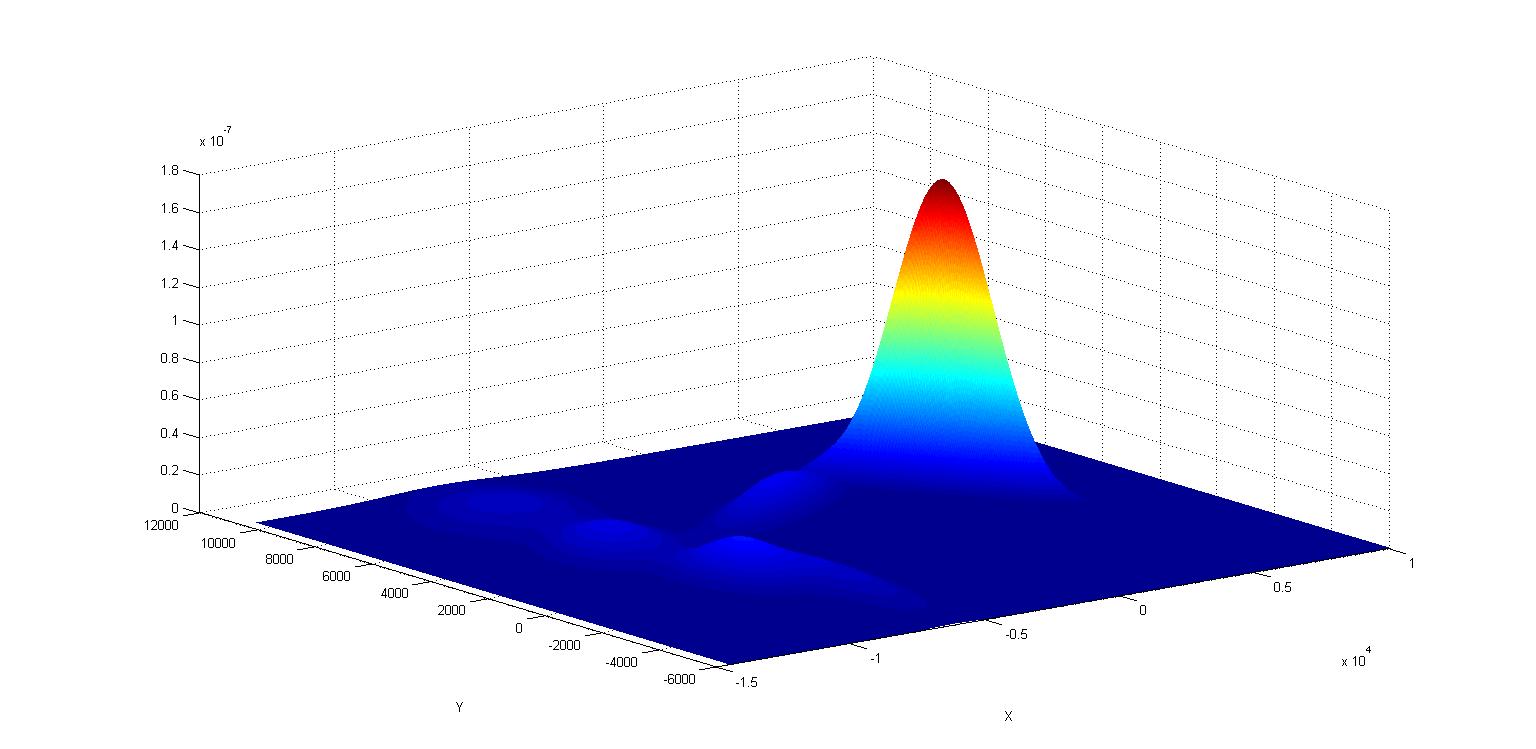}
      \caption{$t=T/2$}
                \label{fig:pendangdisp}
        \end{subfigure}
        
        \begin{subfigure}[b]{0.2\textwidth}
                \includegraphics[width=\textwidth]{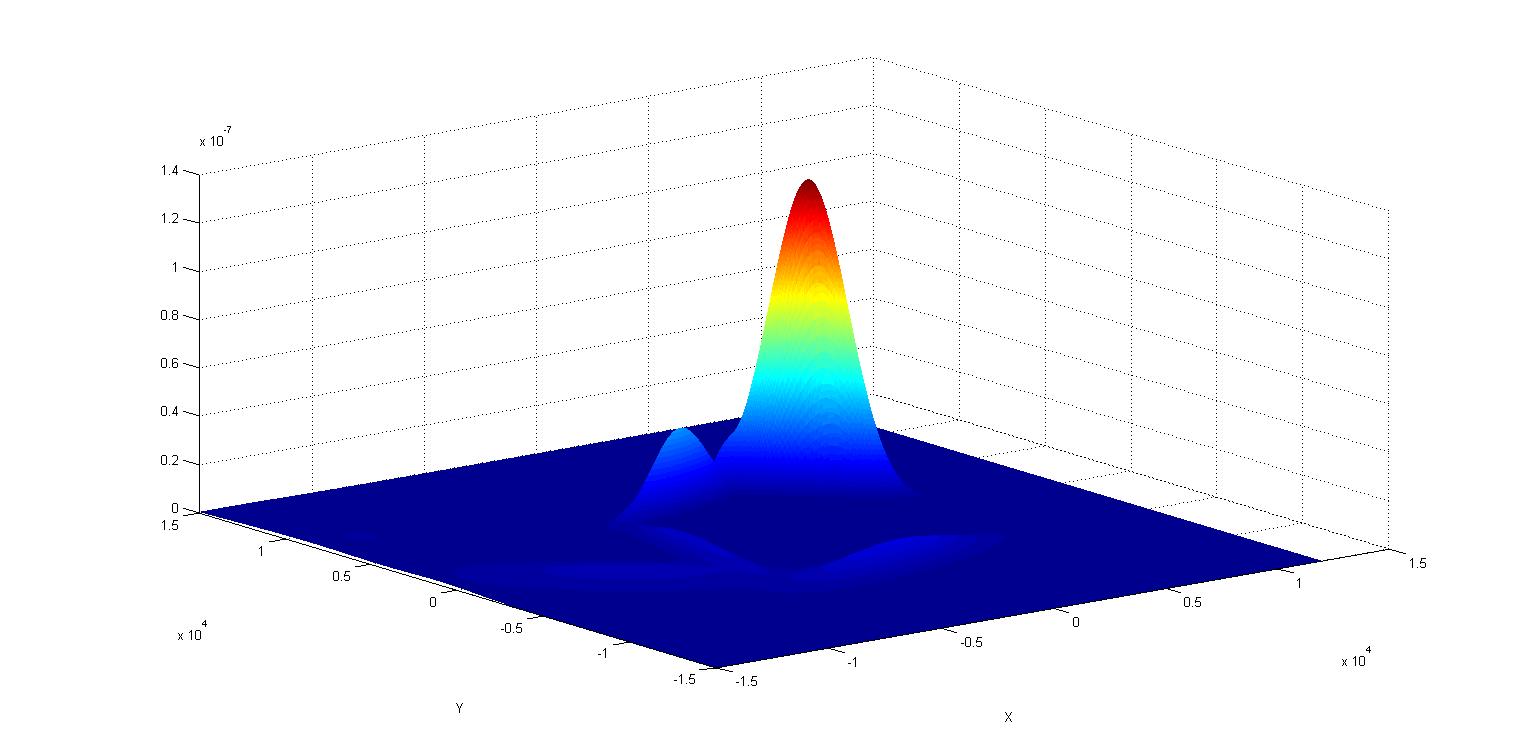}
                \caption{$t=3T/4$}
                \label{fig:pendvel}
                \end{subfigure}
         \begin{subfigure}[b]{0.2\textwidth}
        \includegraphics[width=\textwidth]{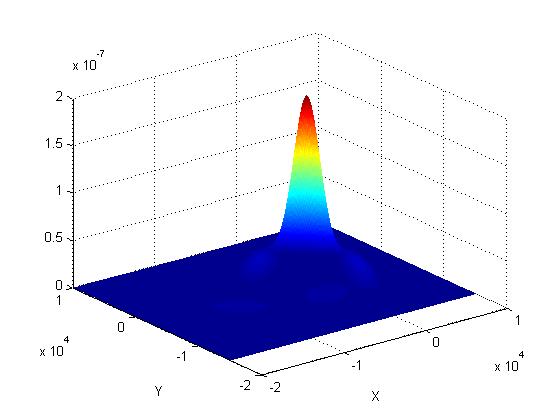}
        \caption{$t=T$}
        \label{fig:pendangvel}
        \end{subfigure}
              \caption{GMM representation of the particle ensemble for times $t>0$ with large initial uncertainty in $V$}
        \label{fig:amocodia}
\end{figure}

 every initial conditions the same as that in the case plotted in figure~\ref{fig:ensemble1} except for the variance in speed which was reduced to 10 m/s.  The distribution of these particles as they orbit the earth is plotted in fig~\ref{fig:margdist}. Since the mechanical energy of the space object increases as the square of the velocity, a larger uncertainty in velocity leads to to a larger variation in the mechanical energy which manifests as a larger uncertainty in the semi-major axis. Thus, a lower initial  uncertainty in velocity may be a great advantage in estimation since the trajectories are closer off and the resulting pdf is less stretched and distorted and more akin to a Gaussian for a longer time period as seen in figure \ref{fig:margdist}. The results obtained by clustering the samples with lower uncertainty in velocity are provided in Table~\ref{tab:lownoise}. Clearly , the number  of Gaussian components and their covariances are smaller  in comparison to the previous case. \\ 
 
\item Negligible process noise: The orbital perturbations that influence the dynamics of space objects are well studied. Consequently, accurate models describing the perturbing forces are  available. Accordingly, the uncertainty involved in the governing equations, i.e., the process noise, can assumed to be negligible when a dynamic model of appropriate accuracy is employed.\\
%

\item Periodicity of the orbits: Orbits occupied by a large class of space objects are periodic or nearly periodic. This implies that their trajectories in the state space approximately retrace their paths when propagated over multiple time periods.
\end{enumerate}

\section{Filter Design}

 A space object tracking scenario  usually involves  one or more  sensors that are employed to record observations. Since the sensing resources  have only a  limited field of view (FOV), only  part of the sky is visible at any instant. Thus, an orbiting object could be outside the field of view
\begin{table}[H]
\centering
\caption{}
\label{tab:lownoise}
   \begin{tabular}{ |c|c|l| }
     \hline
  \multicolumn{3}{|c|}{Low noise} \\
  \hline
  time & No. of modes & \multicolumn{1}{c|}{trace ($\times 10^2$)} \\
\hline
  0 & 1 & $3.0938 \times 10^-2$ \\
 $1500$ & $1$ & $2.956 \times 10^{-1} $ \\
  $3000$ & $2$ & $1.856$ \space{} $1.823$  \\
 $4500$ & $3$ & $3.850$ \space{} $3.593$ \space{} $3.941$\\
$6000$ & $2$ & $2.659$ \space{} $2.903$\\  \hline
\end{tabular}
\end{table}  
  for long times during which no measurements would be available. Nevertheless, given a sufficiently high probability of detection, it would be safe to assume that the sensor would generate frequent measurements of the object while its inside the field of view. Under these circumstances, a UKF-PF hybrid tracking scheme may be conceived to harness the advantages of both filters. During its flight inside the field of view,  the UKF can be employed to estimate the state of the space object. The between-measurement distortion of the state pdf from Gaussian will be limited if frequent measurements are available. It is also computationally efficient as the unscented transform requires only 2n+1 sample points for an n dimensional estimation problem. Moreover, state estimation with UKF is free from resampling procedures that are customary to sequential Monte Carlo methods.\\

  Once the space object is outside the field of view, the pdf is extensively distorted since no further measurements will be recorded till the object re-enters the FOV. As it is free from restrictive assumptions, it is more advantageous to use a particle filter here especially since the actual pdf may be captured to any degree of accuracy by increasing the number of particles. Once the object exits the FOV and the particles are sampled, they are propagated through time while keeping the individual weights $ W(x^{i})$ constant until the next measurement is recorded, i.e. when the object re-enters the FOV. The negligible process noise and near periodic orbits help enable the propagation of the state pdf for extended time periods when the measurements are sparse without risking particle depletion which is a problem often observed in typical PF implementations\cite{Aru}\cite{Docet}. As objects are first detected when they are within the FOV, the uncertainty in velocity will be reduced before the state 
  pdf is propagated through the particle filter.\\
 
 Employing a hybrid UKF-PF filter warrants laying out the procedure for transition from UKF to PF and vice versa. To transition from UKF to PF, the necessary number of particles can be directly sampled from the prior pdf given by the UKF. This has to be carried out once the estimated position of the object is outside the field of view. Transitioning from PF to UKF on the other hand requires a Gaussian pdf to be retrieved back from the ensemble of particles. Once the object re-enters the FOV and the first measurement is registered, the particles are assigned weights based on their respective likelihood derived from the measurement model.

 The likelihoods for each particle $X_{i}$ at instant $k$ may be computed as
\begin{equation}
W(X_{i})=P_{\nu}(Z_{k}-H(X_{i})),
\end{equation}
where $Z_{k}$ is the measurement recorded at that instant. The mean and covariance of this weighted sample can be used to parameterize the necessary Gaussian pdf which may  be computed as
\begin{equation}
\mu=\Sigma_{i=1}^{N} X_{i}W(X_{i}),
\end{equation}

\begin{equation}
C=\Sigma_{i=1}^{N}\dfrac{ W(X_{i})(X_{i}-\mu)^{2}}{1-\Sigma_{j=1}^{N}W(X_{i})^{2}}.
\end{equation}

\begin{figure}[H]
\includegraphics[width=0.5\textwidth,height=0.15\textheight]{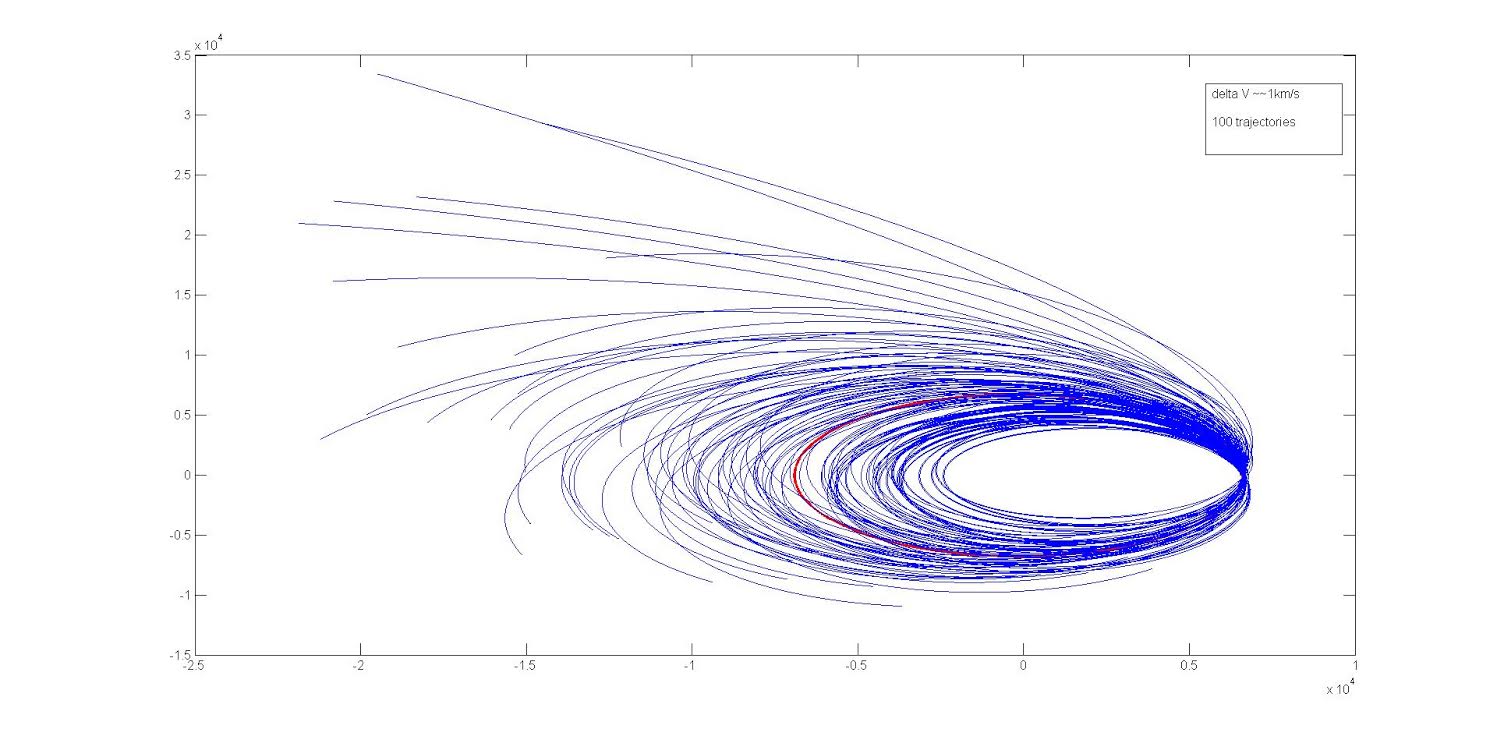}  
\caption{Ensemble of trajectories with large velocity uncertainty}
\label{fig:trajectories}
\end{figure}

  While there is a chance that a significant fraction of the particles may undergo weight degeneration, a  lower initial uncertainty in velocity keeps the trajectories closer and the risk of weight depletion can be avoided. This is also the case that is encountered in most practical situations. For instance, when propagated with an initial uncertainty of 10m/s in each direction, the particles did not suffer from weight depletion even after significant time spans. It should be emphasized that if the process noise had been large, the particles would have been considerably scattered irrespective of the magnitude of initial uncertainty. Hence, the absence of  model uncertainty plays an important role in preventing particle depletion. Fig~\ref{fig:trajectories} shows the trajectories followed by particles sampled from an initial Gaussian distribution with mean $[6600cos\pi /12 \hspace{6pt} 0 \hspace{6pt} 6600sin\pi/12 \hspace{6pt} 0 \hspace{6pt} 7.8848\hspace{6pt} 0]^{T}$  for the 2D Keplerian problem. Even though the particles appear to diverge in the beginning, they come back to their respective initial states on completion of a full time period.
Clearly the trajectories appear denser near the initial state and are less likely to encounter weight depletion at nearby regions.Hence the periodicity of the orbit also plays a significant role in preventing weight depletion. 
 Based on this observation, the following proposition regarding particle depletion in such noise free periodic systems may be made.
\begin{proposition}
 Consider a  periodic dynamical system governed by 
\begin{equation}
\dot{X}=f(X),
\end{equation} 
with an uncertain initial state characterized by a Gaussian pdf with mean and covariance ($S_{0},P$). Given a measurement model  
\begin{equation}
Y=g(X)+\nu,
\end{equation}
associated with this system, where $\nu$ is a zero mean Gaussian noise term with covariance R, a lower bound on the probability that the measurement likelihood of any particle that is sampled from the initial Gaussian pdf, would be above a given threshold $b$, after a full period $T(S_{0})$, of the mean trajectory, is given by the total probability enclosed inside an $m-\sigma$ contour  ellipse of a Gaussian pdf parameterized by the covariance matrix C given by $C=2MPM^{T}+R$. Here $M=\dfrac{dg}{dX}\vert_{X=S_{0}}(I-f(So).\dfrac{dT}{dS})$, where $S_{o}$ stands for the mean initial state $S$, $ T(S)$ for period of the orbit as a 
function of initial state $S$ and  m may be computed as  $m=\sqrt{\frac{\alpha_{min}}{\lambda_{max}}log(\frac{1}{b\sqrt{2 \pi det(R)^{n}}})}$, where $\alpha_{min}$ is the smallest eigenvalue of $R$ and $\lambda_{max} $ is the largest eigenvalue of C. 
\end{proposition}
Clearly, a lower value of $m$ would indicate a higher risk of weight degeneration. From the given expression, it can be observed that $m$ is directly proportional to the smallest eigenvalue of R, $\alpha_{min}$, and inversely proportional to the largest eigenvalue of the matrix M, $\lambda_{max}$. Hence, a smaller measurement covariance and a larger initial state uncertainty would increase the risk of particle depletion. The matrix $M$ indicates the transformation that the covariance matrix P undergoes over a full period of the mean trajectory.   From equation of the matrix M, it can be concluded that the sensitivity of time period w.r.t. the initial condition, $\dfrac{dT}{dS}$, plays an important role in determining the risk of weight degeneration. For  gravitational systems, the square of time period is proportional to the cube of the semi major axis which in turn is a function of the total energy. But since mechanical energy varies with the square of the velocity, it can be concluded that the time period of the gravitational system is sensitive towards the uncertainty in velocity.  Hence, shrinking the uncertainty in velocity would keep the periods of the sampled trajectories closer off. As a result, particles along these trajectories would traverse through their initial states at almost the same times, reducing the risk of weight degeneration at that point.\\

\indent The above result is arrived at based on the fact that variation in the relative position of two points with similar periods over a full time period of either of these particles  can be computed as a function of the difference in their orbital periods. This assumes that the time period $T$ is a differentiable function of the initial state $S$. From this, the set of points for which the given condition on the likelihood function would hold, after a full time period, of the mean trajectory, may be computed. Then the lower bound that is presented in the given proposition may be computed as the probability of sampling particles from this particular set given the initial Gaussian pdf $\mathcal{N}(S_{0},P)$. A derivation of the above result is presented in the appendix \\

 \indent Fig~\ref{fig:transition} shows the particle distribution in the x-y plane during the transition from PF to UKF. As expected, the particle distribution prior to registering the measurement appear stretched and distorted. Once the weights are computed, the particles are resampled to generate a new ensemble. If the measurement uncertainty is small compared to the covariance of the state pdf, the resampled set would be less stretched and distorted as it can be seen from fig~\ref{fig:transition}. Hence, it is reasonable to approximate the resampled pdf with a single Gaussian component. Based on this observation,the mean and covariance of this new pdf is employed in the UKF based estimation process that follows the transition. The steps involved in implementing the UKF-PF hybrid filter are described in algorithm~\ref{algo}.
 \begin{figure}[H]
\includegraphics[width=0.45\textwidth,height=0.15\textheight]{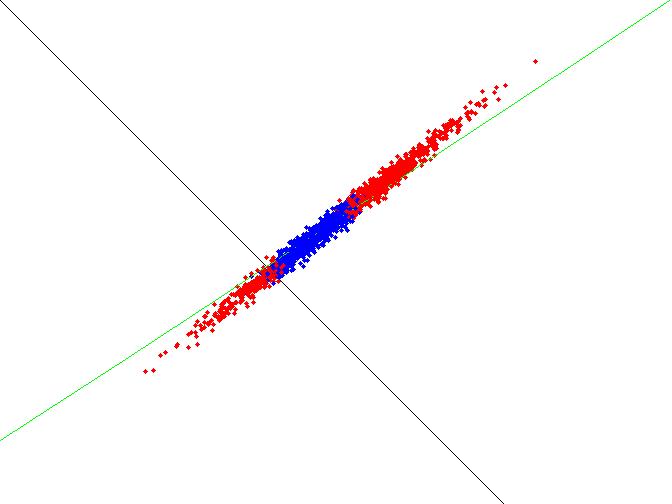} 
\caption{Ensemble of particles before(red) and after(blue) resampling}
\label{fig:transition}
\end{figure}

\section*{Simulations and Results}

The nonlinear filter developed in the previous section is employed to estimate the state of a space object in the low earth orbit (LEO). The resultant  acceleration experienced by an object in LEO may be computed as 
\begin{equation}
\ddot{r}=a_{g}+a_{j_{2}}+a_{d}.
\end{equation}
The state of the system is taken to be
\begin{equation}
X=[x_{1}\hspace{3pt}x_{2}\hspace{3pt} x_{3}\hspace{3pt} \dot{x_{1}}\hspace{3pt} \dot{x_{2}}\hspace{3pt} \dot{x_{3}}]^{T}.
\end{equation}
Any model uncertainty present is quantified with a process noise term defined by zero mean Gaussian pdf with a covariance $10^{-10}\mathbf{I_{6}}$. The process noise term is modelled  as an acceleration resulting from an unknown forcing and hence added in the time update equations for the system state. 
Performance of the proposed filtering scheme  is assessed by simulating the orbital dynamics for 10 periods and evaluating  the root mean squared error (RMSE) matrix of the estimated state and comparing with posterior Cramer-Rao lower bound(PCRB)\cite{Tich}. The PCRB establishes a lower bound on the mean square error for the filtering problem.
As a result, the matrix difference, A of the root mean square error(RMSE) matrix and PCRB matrix is always positive semi-definite, i.e, zero is a lower bound for the eigenvalues of A. The spectral norm of a matrix A is defined by
\begin{equation}
\parallel A \parallel_{2}=(\lambda_{max}(A^{H}A))^{1/2}.
\end{equation}
Hence, a lower spectral norm suggests a better filtering performance. Also, the smallest eigenvalue of the mean
\begin{algorithm}[H]
\caption{UKF-PF Hybrid filter for space object tracking}

$S_{1}$ : PDF in functional form,
$S_{2}$ : PDF as ensemble \\
$C(X)$ : Boundary of FOV\\
$\eta_{0}$ : Probability of detection

Initialize: $P(X)=P_{0}(X), S=S_{0}$\\
At $t_{k}$
\begin{algorithmic}[1]
\IF {$S=S_{1}$}
\IF{$C(x_{k})<= 0$ }
\STATE{Use UKF}
\STATE{SET $S=S_{1}$}
\ELSE
\STATE{SAMPLE FROM P(X)}
\STATE{USE PF}
\STATE{SET $S=S_{2}$}
\ENDIF
\ELSE
\IF{$C(x_{k})<= 0$ \AND $\eta > \eta_{0}$ }
\STATE{USE PF}
\STATE{COMPUTE PARTICLE WEIGHTS}
\STATE{COMPUTE P(X)}
\STATE{SET $S=S_{1}$}
\ENDIF
\ELSE
\STATE{USE PF}
\STATE{SET $S=S_{2}$}
\ENDIF
\end{algorithmic}
\label{algo}
\end{algorithm}

 square error matrix, $\lambda_{min}$, should always be non-negative. An exception would indicate that the mean square error in the state estimate is underestimated.
Also employed is the normalized estimation error squared (NEES) test for evaluating filter consistency\cite{Bail}. This involves the calculation of  a quantity $\beta_{k}$, which is  defined as 
\begin{equation}
\beta_{k} = ( x_{k}- \mu_{k|k}) )^T \bf{P}^{-1}_{k|k}( x_{k}- \mu_{k|k}) ).
\end{equation}
 For a six dimensional random variable, expected value of $\beta_{k}$ is 6, while 90 per cent of its  probability is concentrated between the values 1.635 and 12.592. If $\beta$ assumes a value above 12.592, then it is more likely that the covariance $P_{k|k}$ was underestimated, in other words, the estimates were very optimistic. Similarly, if the value of $\beta$ is lower than 1.635, it is likely that the covariance was over estimated, i.e. the estimate was conservative.

The UKF-PF hybrid filter is implemented as the state estimation scheme in the following test cases.\\

\begin{figure}[h]
  \begin{subfigure}[b]{0.45\textwidth}
   \includegraphics[width=\textwidth]{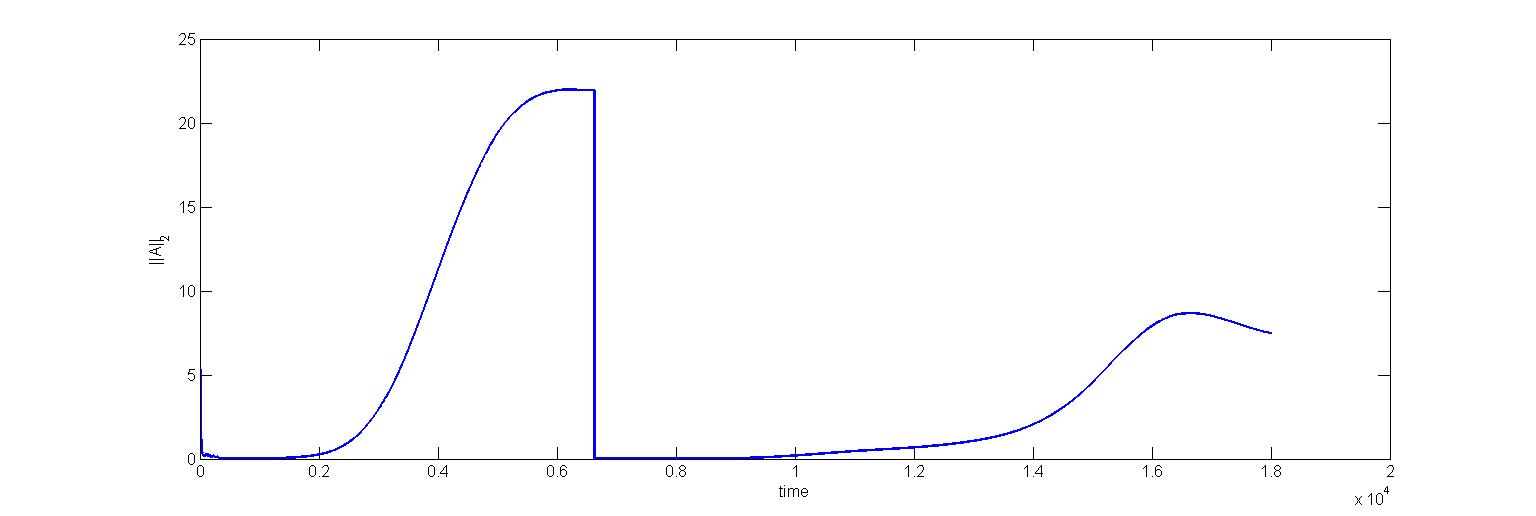}
   \caption{}
 \label{fig:penpio}
 \end{subfigure}
        \begin{subfigure}[b]{0.45\textwidth}
                \includegraphics[width=\textwidth]{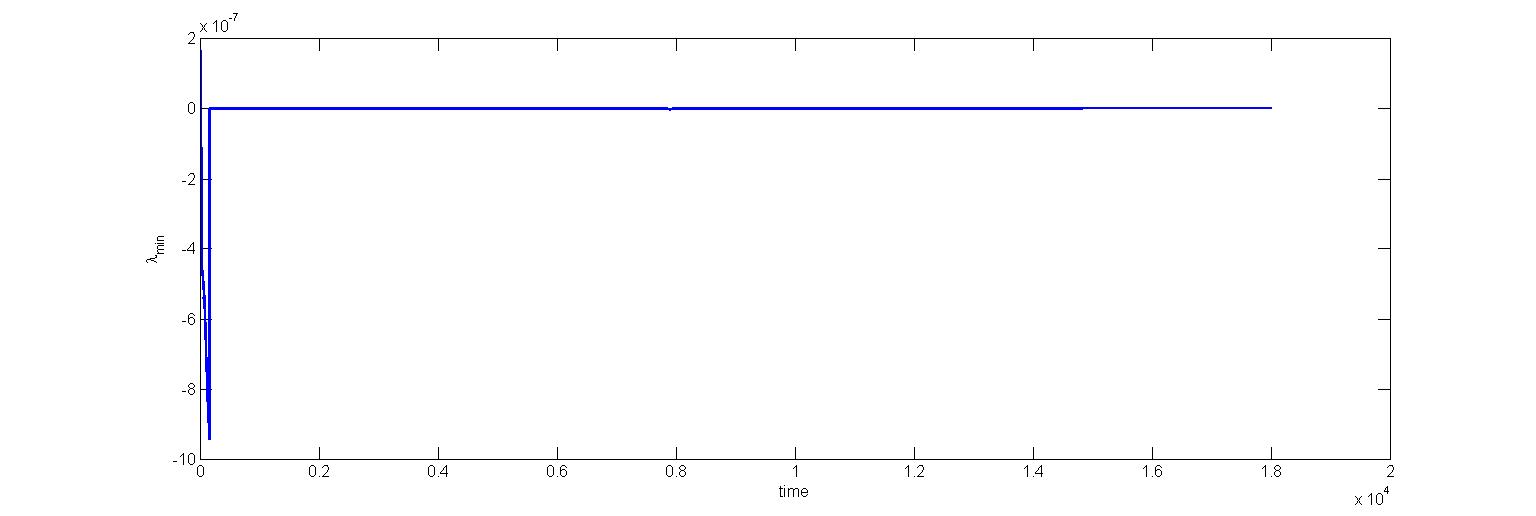}
                 \caption{}               
                \label{fig:pendiyangdisp}
        \end{subfigure}

        \begin{subfigure}[b]{0.45\textwidth}
        \includegraphics[width=\textwidth]{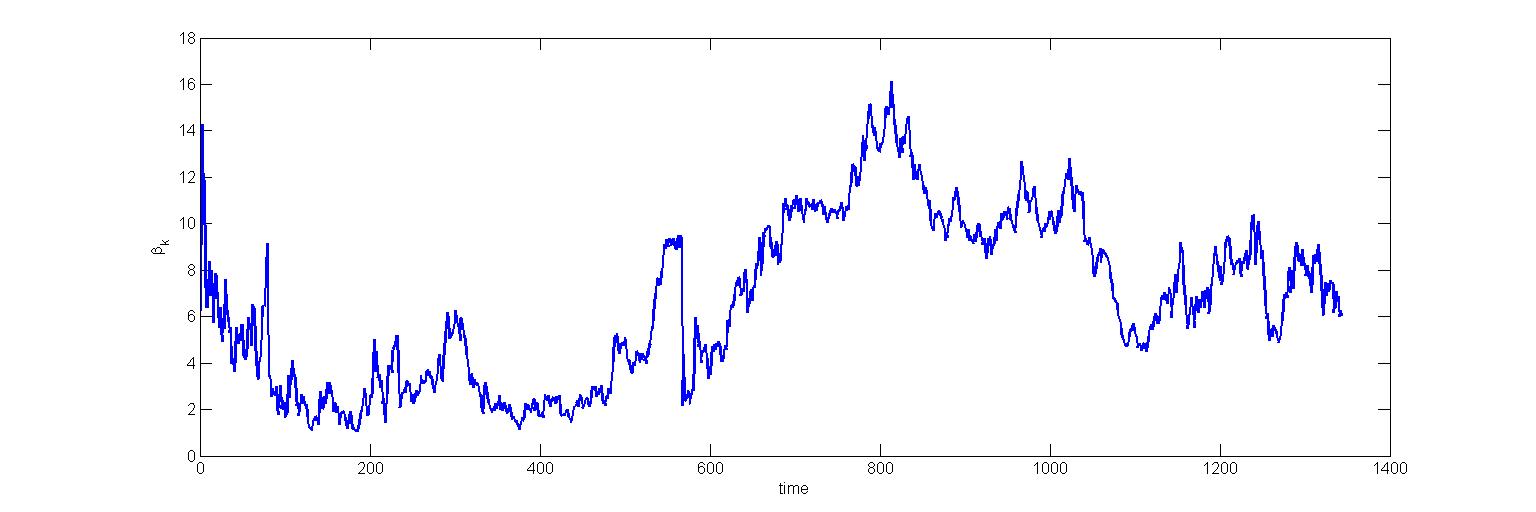}
          \caption{}
          \label{fig:pendvoell}
     \end{subfigure}
         \caption{Test case 1 results: (a) Spectral norm($\vert \vert A \vert \vert_{2}$)from PCRB (b) $\lambda_{min}$ from PCRB (c) NEES plot}
        \label{fig:margdist5}
\end{figure}

\textbf{Case 1}:

In this case, the initial state of the space object is set at\\ 
$S_{0}=\begin{bmatrix}7800 \hspace{7pt} 0\hspace{7pt} 0 \hspace{7pt} 0 \hspace{7pt} 6.8443cos(pi/4) \hspace{7pt} 6.8443sin(\pi/4)\end{bmatrix}^{T}$
 where the lengths and speeds are in km and km/s respectively. This is a $45$ degree inclined low earth orbit with a period 6080 s and eccentricity 0.0833. There is $5\hspace{2pt} km$ standard deviation in the initial position estimate and $1\hspace{2pt} m/s$ standard deviation in velocity estimate in each directions. The orbital dynamics and the estimation scheme was simulated over a time period of 5 hours. The filtering results for the above space object for the said period is studied in terms of PCRB and NEES and plotted in fig~\ref{fig:margdist5}. In fig ~\ref{fig:penpio} spectral norm of the matrix difference between RMSE matrix and PCRB is plotted. The value of spectral norm is close to zero when measurements are available. But, once the space object moves out of the FOV of the sensor, the amount of available information reduces steadily as signified by the upswing in the spectral norm. The lowest eigen value of the matrix A is a non-negative infinitesimal during most of the times. The exception occurs during the beginning of the simulation when both RMSE and PCRB are almost equal\cite{Tich} and their infinitesimal difference is calculated to be a small negative quantity which we attribute to numerical round off error. Figure ~\ref{fig:pendvoell} shows the results of NEES test for test case 1. It should be noted that NEES is evaluated only during those instants when measurement is available. Hence the time axis in \ref{fig:pendvoell} indicates time steps with measurement updates. Evidently, the estimation is not too optimistic or conservative most of the time. Over 1344 times instants at which NEES was computed, only 111 (8.25 per cent) had the value of $\beta_{k}$ beyond the 90 per cent bounds. Even though this is a six dimensional dynamic system, only a small number of particles (2000) were used in this simulation.\\

\begin{figure}[h]
       
        \begin{subfigure}[b]{0.45\textwidth}
                \includegraphics[width=\textwidth]{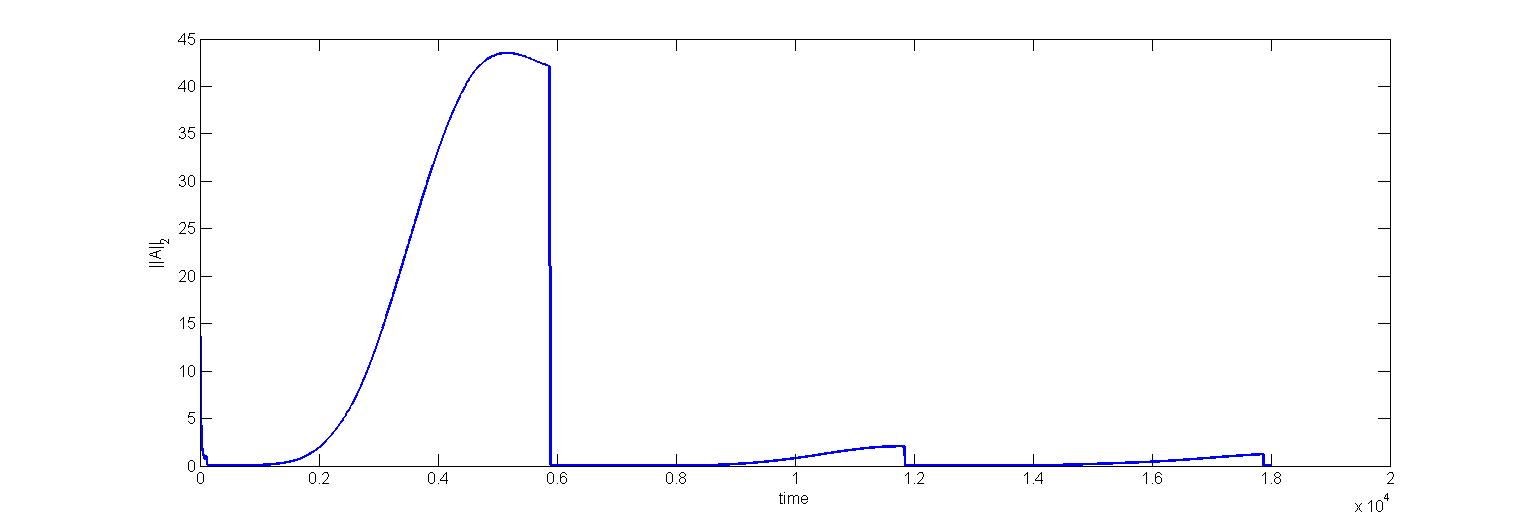}
                
                \label{fig:peyi}
        \end{subfigure}
        \begin{subfigure}[b]{0.45\textwidth}
                \includegraphics[width=\textwidth]{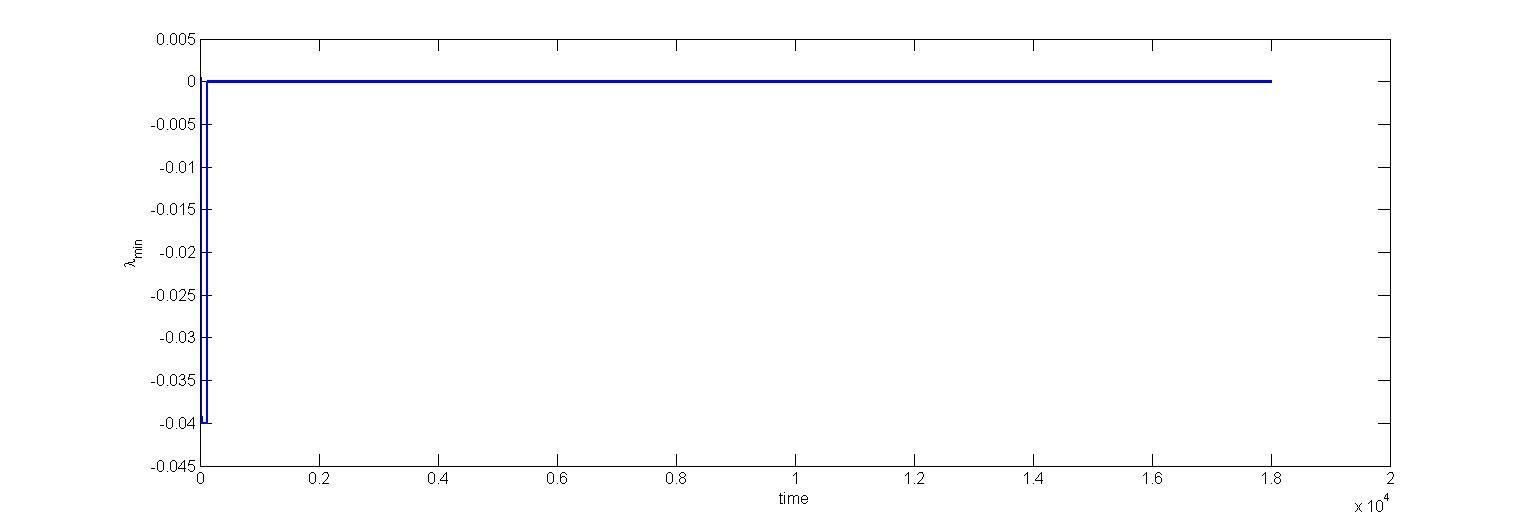}
                                
                \label{fig:pendiyaisp}
        \end{subfigure}
        
        \begin{subfigure}[b]{0.45\textwidth}
                \includegraphics[width=\textwidth]{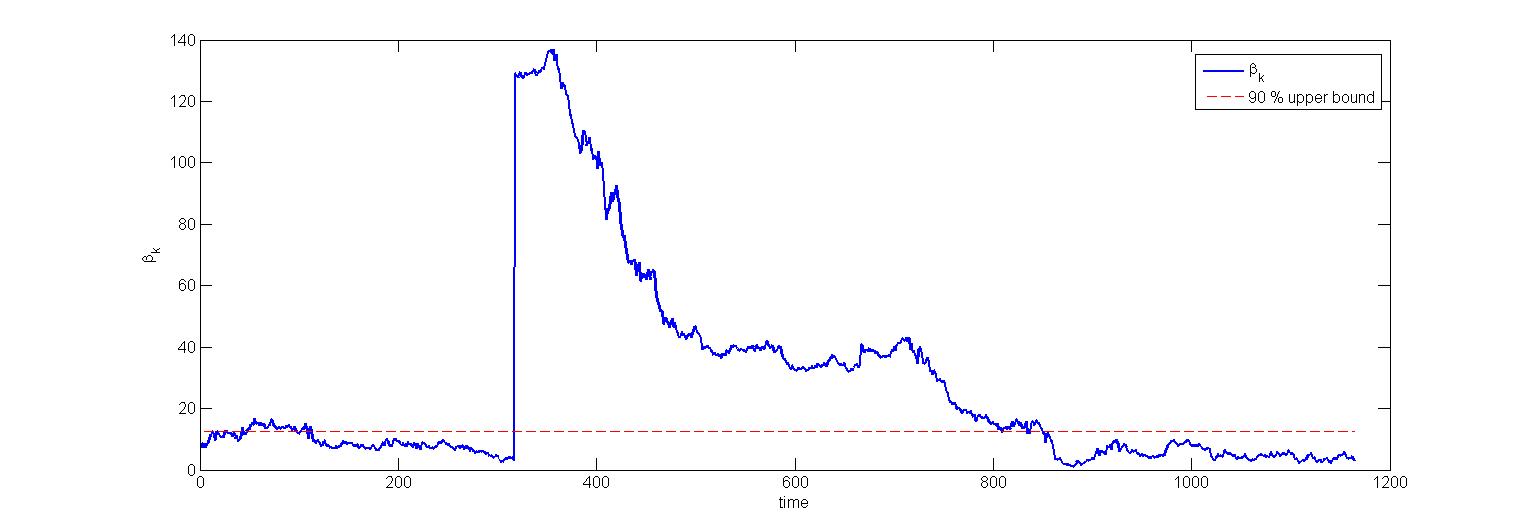}
                
                \label{fig:penoel}
                \end{subfigure}
         \caption{Case 2 results: (a) Spectral norm ($\vert \vert A \vert \vert_{2}$) from PCRB (b) $\lambda_{min}$ from PCRB (c) NEES plot}
        \label{fig:margdist61}
\end{figure}

\textbf{Case 2}:

In this case, the initial state of the space object is set to \\
$S_{0}=\begin{bmatrix}6800 \hspace{7pt} 0\hspace{7pt} 0 \hspace{7pt} 0 \hspace{7pt} 7.5989cos(\pi/30) \hspace{7pt} 7.5989sin(\pi/30)\end{bmatrix}^{T}.$
This is orbit has a lower inclination ($\pi/30$ radian from the equatorial plane) in comparison to case 1 and has a time period of 5580.5 s  The initial uncertainty in position was kept at 2 km but that in velocity was raised to 200 m/s in each directions to study the performance of the hybrid filter in large noise scenarios. The simulation was run over a period of $5$ hours. The result obtained from test case 2 is given in Fig~\ref{fig:margdist61}. The proposed filter offers good estimation performance in terms of the spectral norm of PCRB metric($\vert\vert A\vert\vert_{2}$). The growth in $\vert\vert A\vert\vert_{2}$ during the flight of the object outside FOV is seen to reduce gradually as the total number of recorded measurement grows. From Fig~\ref{fig:penoel}, it can be seen that the value of $\beta$ spikes beyond the 90 percent upper bound suddenly at the point where the first transition from PF to UKF occurs. indicates that the uncertainty in the state estimate is underestimated. Due to the larger initial uncertainty in velocity, the particle ensemble stretches and distorts more in test case 2 when the object is outside FOV. Consequently, assuming a single mode Gaussian pdf for the state pdf at the transition point would result in underestimating the actual uncertainty present in the state estimate. However, the value of $\beta_{k}$ reduces gradually and no further spikes are observed even though there are two more PF to UKF transitions as can be seen the PCRB plots. As more measurements are recorded, the uncertainty in the state pdf reduces and $\beta_{k}$  shrinks and stays mostly within the 90 percent bounds. This indicates that given sufficient number of measurements, the hybrid filter offers reliable and consistent estimation performance even if there is large initial uncertainty

\section*{Conclusions}
 The design and application of a hybrid UKF PF based estimation scheme for tracking space objects has been presented. The dynamics of space objects under the effects of $J_{2}$ perturbation and atmospheric drag is considered. The conventionally employed nonlinear filters such as EKF, UKF etc. require the restrictive assumption that the state pdf of a nonlinear dynamical system remains Gaussian at all times. This is particularly problematic when the pdf is distorted under nonlinear transformations for extended periods of time. \\
\indent The proposed hybrid filter employs the UKF for tracking when the relevant space object is inside the field of view and measurements are recorded. In order to handle the non-linear distortion outside the FOV, the tracking scheme transitions to PF as the object moves  outside FOV.  It was found that a smaller uncertainty in velocity could prevent particle depletion. Also, limited model uncertainty in the orbital problem, and the near periodic dynamics involved, play no less a role in preventing particle depletion. It is also found that while the pdf undergoes extensive nonlinear distortion when there are no measurement updates, a small measurement uncertainty and large probability of detection could facilitate the use of UKF without incurring a large 

error while within the FOV.
The proposed filtering scheme was employed to estimate the state of a space object in inclined low earth orbits and the performance is studied in terms of the PCRB and NEES metrics. It is demonstrated that the hybrid filter is a fast and computationally inexpensive space object tracking scheme and offers excellent performance regardless of the orbital inclination or initial uncertainty as long as sufficient number of measurements are recorded within the FOV. Also the new filter was not observed to suffer from particle depletion in any of the test cases despite using a relatively small number of particles.\\
\section*{Appendix I}
Consider a dynamical system with known governing equation and initial state distribution. Uncertainty in the initial state will propagate through the state space in time, based on the time evolution of the system dynamics. Hence, the state of the system would be unknown $\forall t>0$ and will be specified through a pdf unless new knowledge is incorporated. Uncertainty propagation can be studied by sampling initial states based on the initial state pdf. The trajectories arising from this sample would be studied for $\forall t>0$. In the limit as the number of sampled states, $n_{S}$, reach infinity, the sample would become an exact representation of the actual distribution of trajectories. The sampled particles are updated or re-sampled whenever an observation is recorded. During the re-sampling process, trajectories arising from many samples would be dismissed as improbable, based on the disparity between the actual and calculated (from the observation model) values of observation. Sometimes, the number of trajectories that are retained after the resampling process would be too small that the resulting estimate would underestimate the actual uncertainty involved in the system. This problem is called particle depletion. In the case of periodic dynamic systems with Gaussian uncertainties(for initial state and measurement model) a lower bound for the probability that any sampled trajectory would retain a minimum likelihood after a full time period of the mean trajectory, can be calculated.

Let $f$ be a dynamical system given as
\begin{equation}
\dot{X}=f(X).
\end{equation}
       If the solution of the system is periodic for a given initial condition S, then the resulting trajectory Q, specified as $Q(t,S)$ satisfies the following. 
       \begin{eqnarray}
        Q(t,S)= Q(t+T(S),S),\\
       Q(-t,S)= Q(T(S)-t,S).
       \end{eqnarray}

Let the initial state be specified with a Gaussian pdf with mean $S_{0}$ and covariance $P$. Suppose the measurement model is given as 
\begin{equation} Y=g(X)+\nu,
\end{equation}
where $\nu$ is zero mean Gaussian white noise with a covariance $R$. Consider the mean trajectory starting from $Q(0,S_{0})$ and any randomly sampled trajectory Q(S). Let $T_{0}$ be the time period of the mean trajectory, i.e. $T_{0}=T(S_{0})$. If at $t=0$ the sampled points are at $S_{0} (=Q(0,S_{0}))$ and  $S(=Q(0,S))$, then after a time $T_{0}$, they will be at $Q(T_{0},S_{0})$ and $Q(T_{0},S)$. Since both trajectories are periodic, we get 
\begin{eqnarray}
Q(T_{0},S_{0})=Q(T_{0}-T(S_{0}),S_{0})=Q(0,S_{0}),
\end{eqnarray}
\begin{equation}
Q(T_{0},S)= Q(T_{0}-T(S),S).
\end{equation}

Now let the actual underlying state of the system be $S^{*}$. Then after the time $T_{0}$, the actual state would be $Q(T_{0}-T(S^{*}),S^{*})$. Since the measurements registered originate from the actual state, we get
\begin{equation}
y=g(Q(T_{0}-T(S^{*}),S^{*}))+\nu .
\end{equation}
Once this measurement is registered, the likelihoods for different trajectories are computed as 
\begin{equation}
\label{likely}
L(S)= \frac{1}{\sqrt{2 \pi det(R)^{n}}} e^{(\Delta y)^{T} R^{-1}(\Delta y)},
\end{equation}
where $\Delta y= g(Q(T_{0},S)-y$. Now,
\begin{equation}
\label{long}
\begin{split}
g(Q(T_{0},S)-y  = g(Q(T_{0},S)- g(Q(T_{0}-T(S^{*}),S^{*}))-\nu \\
            =  g(Q(T_{0}-T(S),S))- g(Q(T_{0}-T(S^{*}),S^{*}))-\nu \\
            =  g(Q(T_{0}-T(S),S))- g(Q(0,S_{0}))
               +g(Q(0,S_{0}))\\- g(Q(T_{0}-T(S^{*}),S^{*}))-\nu .
\end{split}
\end{equation}

 Linearizing $g(x)$ at $S_{0}(=Q(0,S_{0}))$ with Taylor series gives 
\begin{equation}
\begin{multlined}
   g(Q(t,S))\\ \approx G(Q(0,S_{0}))+\dfrac{dg}{dX}|_{X=S_{0}}(\dfrac{dQ}{dt}|_{X=S_{0}}.t+\dfrac{dQ}{dS}|_{X=S_{0},t=0}\Delta S)\\
            = G(Q(0,S_{0}))+\dfrac{dg}{dX}|_{X=S_{0}}(f(So).t+I\Delta S).
\end{multlined}             
\end{equation}

From this, we find
\begin{equation}
\begin{multlined}
g(Q(T_{0}-T(S),S))- g(Q(0,S_{0}))\\=\dfrac{dg}{dX}(f(So).(T_{0}-T(S))+I\Delta S).
\end{multlined}
\end{equation}
Since time period $T$ is a function of the initial state $S$, it can be linearized to find 
\begin{equation}
T_{0}-T(S)=\dfrac{dT}{dS}|_{X=S_{0}}.\Delta S.
\end{equation}
Using this, the eqn~\ref{long} can be reduced to
\begin{equation}
\begin{split}
g(Q(T_{0},S)-y =\dfrac{dg}{dX}|_{X=S_{0}}(I-f(So).\dfrac{dT}{dS}|_{X=S_{0}})\Delta S_{1}\\ - \dfrac{dg}{dX}|_{X=S_{0}}(I-f(So).\dfrac{dT}{dS}|_{X=S_{0}})\Delta S_{2}-\nu .
\end{split}
\end{equation}
where $\Delta S_{1}=S-S_{0}$ and $\Delta S_{2}=S^{*}-S_{0}$. Hence
\begin{equation}
\label{long1}
\begin{split}
g(Q(T_{0},S)-y=\dfrac{dg}{dX}(I-f(So).\dfrac{dT}{dS})(\Delta S_{1}-\Delta S_{2})-\nu .
\end{split}
\end{equation}
all functions evaluated at $X=S_{0}$. For a given $S_{0}$, the coefficient of $\Delta S_{1}-\Delta S_{2}$ in eqn~\ref{long1} is a constant. Hence it can be written that 
\begin{equation}
\label{long2}
\begin{split}
\Delta y=M(\Delta S_{1}-\Delta S_{2})-\nu,
\end{split}
\end{equation}
where $M=\dfrac{dg}{dX}(I-f(So).\dfrac{dT}{dS})$.
 \\Now, consider the likelihood function discussed in eqn~\ref{likely}. For the likelihood function $L$ to be bounded by a lower bound $b$
\begin{equation}
\label{inequal1}
\begin{split}
\frac{1}{\sqrt{(2 \pi)^{k} |R|}} e^{-[(\Delta y)^{T} R^{-1}(\Delta y)]} > b,\\
(\Delta y)^{T} R^{-1}(\Delta y)<-log(b\sqrt{(2 \pi)^{k} |R|}),\\
(\Delta y)^{T} R^{-1}(\Delta y)<log(\frac{1}{b\sqrt{(2 \pi)^{k} |R|}}).
\end{split}
\end{equation}
The term on the LHS in the last line of the the inequality~\ref{inequal1} is a quadratic form in the random vector $g(Q(T_{0},S)-y$. Since $R^{-1}$ is the inverse of the measurement noise covariance all its eigenvalues are positive real numbers. Hence, from the theory of quadratic forms, the solution set of this quadratic form is an ellipsoid. The principal axes of the ellipsoid is given by the square root of the eigenvalues of R. When the value of this quadratic form equals the term on RHS,
\begin{equation}
\Delta y^{T}R^{-1}\Delta y=log(\frac{1}{b\sqrt{(2 \pi)^{k} |R|}}) .
\end{equation}
By diagonalizing the covariance matrix R, the $n-\sigma$ ellipsoid for which this relationship holds can be calculated and it is not hard to find that 
\begin{equation}
\label{neq}
\begin{split}
n^2=log(\frac{1}{b\sqrt{(2 \pi)^{k} |R|}}),\\
n=\sqrt{log(\frac{1}{b\sqrt{(2 \pi)^{k} |R|}})}.
\end{split}
\end{equation}
Thus for all $\Delta y$ inside this $n-\sigma$ ellipse the inequality in \ref{inequal1} would hold.
 The random variable $\Delta S_{2}$ in eqn~\ref{long2} is normally distributed with zero mean and a covariance $P$. Also the noise term $\nu$ is disributed according to a zero mean Gaussian with a covariance R. Given a particular initial condition $S$ value of $M\Delta S_{1}$ is fixed. Thus given a particular $M\Delta S_{1}$ mean of the random variable $\Delta y$ would be 
 \begin{equation}
 E[\Delta y]=M\Delta S_{1}.
 \end{equation}
The noise term $\nu$ is independent of the uncertainty in the initial state. Hence , the covariance of $\Delta y$ given  is given as 
\begin{equation}
C[\Delta y]=MPM^{T}+R.
\end{equation}
Since both $M\Delta S_{2}$ and $\nu$ are Gaussian random variables, a linear combination of the two will also be normally distributed. Let $C_{1}=MPM^{T}+R$. Thus, the conditional pdf of $\Delta y$ would be given as
\begin{equation}
P(\Delta y |V_1)=\frac{1}{\sqrt{(2 \pi)^{k} |C_{1}|}} e^{-[(\Delta y-V_1)^{T} (C_{1})^{-1}(\Delta y-V_1)]},
\end{equation}  
where $V_1=M\Delta S_1$. Now, given that a  particular initial state is sampled ( i.e for a particular value of $M\Delta S_{1}$), the probability that the sampled trajectory would retain the lower bound on the likelihood is equal to the probability that the resulting $\Delta y$ would lie inside the $n-\sigma$ ellipse given in eqn~\ref{neq}. Hence
\begin{equation}
P(\Delta y\in A_0|V_1)=\int_{A_0} P(\Delta y|V_1).d(\Delta y),
\end{equation}

where $A_0$ denotes the region encompassed by the $n-\sigma$ ellipse. The initial state $S$ is randomly sampled from a Gaussian pdf with mean zero and covariance $P$.Hence the mean and covariance of the random variable $V_{1}$ are zero and $MPM^{T}$ respectively. Let $C_{2}=MPM^T$. Hence the pdf of $V_{1}$ is given by
\begin{equation}
P(V_{1})=\frac{1}{\sqrt{(2 \pi)^{k} |C_{2}|}} e^{-[(V_{1})^{T} (C_{2})^{-1}(V_{1})]}.
\end{equation}
Thus the probability that any randomly initial state on propagating would retain the lower bound on the likelihood after a period is given by
\begin{equation}
\label{exp}
\begin{multlined}
P(\Delta y\in A_0)=\int_{\Re^{k}}P(\Delta y\in A_0|V_{1})P(V_{1}).d(V_{1}),\\
           =\frac{\int_{\Re^{k}}\int_{A_0} f(\Delta y, C_1,C_2,V_1).d(V_{1}).d(\Delta y)}{(2 \pi)^{k} \sqrt{|C_{1}||C_{2}|}},
\end{multlined}
\end{equation}
where,
\begin{equation}
\begin{multlined}
f_1(\Delta y, C_1,C_2,V_1)\\ =e^{-[(\Delta y-V_{1})^{T} (C_{1})^{-1}(\Delta y-V_{1})
 +(V_{1})^{T} (C_{2})^{-1}(V_{1})]}.
 \end{multlined}
 \end{equation}
 
 Consider now the term in the exponent. In reference [1] it is shown that 
\begin{equation}
\begin{multlined}
(x-\mu_{1})^{T}(\Sigma_{1})^{-1}(x-\mu_{1})+(x-\mu_{2})^{T}(\Sigma_{2})^{-1}(x-\mu_{2})\\=(x-m)^{T}S^{-1}(x-m)+(\mu_{1}-\mu_{2})^{T}M^{-1}(\mu_{1}-\mu_{2}),
\end{multlined}
\end{equation}
where
\begin{flalign} 
& S^{-1}=\Sigma_{1}^{-1}+\Sigma_{2}^{-1}, \\
& m=S(\Sigma_{1}^{-1}\mu_{1}+\Sigma_{2}^{-1}\mu_{2}), \\
& M=\Sigma_{1}+\Sigma_{2}.
\end{flalign}
 Comparing the terms, we get
 \begin{equation}
 x=V_{1},
  \mu_{1}=0,
  \mu_{2}=\Delta y_{1},
  \Sigma_{1}=C_{2},
  \Sigma_{2}=C_{1}.
 \end{equation}
 Hence
 \begin{flalign}
 & S^{-1} = C_{1}^{-1}+C_{2}^{-1} \\
 & m = S(C_{2}^{-1}.0+C_{1}^{-1}.\Delta y_{1})=SC_{1}^{-1}\Delta y\\
 & M =C_{1}+C_{2}. 
 \end{flalign}
Using this, it can be written that
\begin{figure}[h]
\[\begin{split}
(V_{1})^{T}C_{2}^{-1}(V_{1}) +(\Delta y-V_{1})^{T}(C_{1})^{-1}(\Delta y_{1}-V_{1})\\=
(V_{1}-m)^{T}S^{-1}(V_{1}-m)+(\Delta y)^{T}(C_{1}+C_{2})^{-1})(\Delta y).
\end{split}\]
\end{figure}
 Substituting for the term in the exponent in eqn~\ref{exp} and changing the order of integration, we arrive at
 \begin{equation}
 \begin{multlined}
  P(\Delta y\in A_0)\\
 =\dfrac{\int_{A_0}\int_{\Re^{k}}f_2(\Delta y, C_1, C_2)f_3(V_1,S,m)d(V_{1})d(\Delta y)}{(2 \pi)^{k} \sqrt{|C_{1}||C_{2}|}},
 \end{multlined}
 \end{equation}
 where,
 \begin{eqnarray}
 f_2(\Delta y, C_1, C_2)= e^{-(\Delta y)^{T}(C_{1}+C_{2})^{-1}(\Delta y)},\\
  f_3(V_1,S,m)= e^{(V_{1}-m)^{T}S^{-1}(V_{1}-m)}.
 \end{eqnarray}
 
 Integrating $M\Delta S_{1}$ over $\Re^{k}$ while keeping $\Delta y$ constant, we arrive at
 \begin{equation}
 \label{final}
 \begin{multlined}
 P(\Delta y\in A_0)\\=\dfrac{\sqrt{(2\pi )^{k}|S|}}{(2 \pi)^{k} \sqrt{|C_{1}||C_{2}|}}\int_{A_0}f_2(\Delta y, C_1, C_2)d(\Delta y),\\
 =\sqrt{\dfrac{|C_{1}+C_{2}|}{|C_{1}||S^{-1}||C_{2}|}}\int_{A_0}\dfrac{f_2(\Delta y, C_1, C_2)d(\Delta y)}{\sqrt{(2\pi)^{k}|C_{1}+C_{2}|}},\\
 =\sqrt{\dfrac{|C_{1}+C_{2}|}{|C_{1}||C_{1}^{-1}+C_{2}^{-1}||C_{2}|}}\int_{A_0}\dfrac{f_2(\Delta y, C_1, C_2)d(\Delta y)}{\sqrt{(2\pi)^{k}|C_{1}+C_{2}|}},\\
 =\int_{A_0}\dfrac{f_2(\Delta y, C_1, C_2)d(\Delta y)}{\sqrt{(2\pi)^{k}|C_{1}+C_{2}|}}.
 \end{multlined}
 \end{equation}
 
 By inspecting eqn~\ref{final}, it can be seen that the probability that $\Delta y$ would be enclosed within the necessary $n-\sigma$ ellipse is given by the final Gaussian integral that has zero mean and $C_{1}+C_{2}=2MPM^{T}+R$ covariance. Indeed, it can be concluded that the random variable $\Delta y$ is normally distributed with these parameters. The ellipsoids generated from the covariance of the random variable $\Delta y$ will be centred at the origin and their principal axes will be proportional to the square root of the eigenvalues of the covariance matrix. If the length of the largest principal axis of the $m-\sigma$ ellipse resulting from this covariance is smaller than the shortest principal axis of the $n-\sigma$ ellipse discussed earlier, then all points inside this $m-\sigma$ ellipse would satisfy the inequality given in \ref{inequal1}. Thus 

\begin{equation}
\label{meq}
m=\sqrt{\frac{\alpha_{min}}{\lambda_{max}}log(\frac{1}{b\sqrt{2 \pi det(R)^{n}}})},
\end{equation}
where $\alpha_{min}$ is the smallest eigenvalue of the noise covariance matrix R whereas $\lambda_{max}$ is the largest eigenvalue of the covariance matrix of the random vector $\Delta y=g(Q(T_{0},S)-y$. For any point that lies inside this $m-\sigma$ ellipse, the inequality bounding the likelihood will hold. Note that there may be points outside this $m-\sigma$ ellipse for which the inequality may hold. By taking the largest principal axis of the $m-\sigma$ ellipse smaller than the shortest principal axis of the $n-\sigma$ ellipse, we have achieved a conservative estimate of the set of of points for which the inequality \ref{inequal1} holds. 
\section{ACKNOWLEDGMENTS}
\label{.7}
This work is funded by AFOSR grant number: FA9550-13-1-0074 under the Dynamic Data Driven Application Systems (DDDAS) program.


\addtolength{\textheight}{-12cm}

\makeatother


\end{document}